\documentclass[11pt]{scrartcl}                     % onecolumn (standard format)

\usepackage[english]{babel}
\usepackage[ansinew]{inputenc}
\usepackage{amsfonts}
\usepackage{amsmath}
\usepackage{amsthm}
\usepackage{subfigure} 
\usepackage{graphicx}
\usepackage{psfrag}

\newcommand{\R}{\mathbb{R}}
\newcommand{\N}{\mathbb{N}}
\newcommand{\Z}{\mathbb{Z}}
\newcommand{\im}{\mathrm{im}}
\newcommand{\id}{\mathrm{id}}
\newcommand{\wind}{\mathrm{wind}}
\newcommand{\G}{\mathcal{G}}
\renewcommand{\S}{\mathbb{S}}
\renewcommand{\phi}{\varphi}

\newtheorem{defs}{Definition}[subsection]

\newtheorem{cor}[defs]{Corollary}

\newtheorem{lemma}[defs]{Lemma}

\newtheorem{prop}[defs]{Proposition}

\newtheorem{theorem}[defs]{Theorem}

\newtheorem{example}[defs]{Example}
\newtheorem{bem}[defs]{Remark}

\begin{document}

\title{Extending Immersions into the Sphere}
\author{Dennis Frisch}

\maketitle

\begin{abstract}
We study the problem to extend an immersed circle $f$ in the sphere $\S^2$ to an immersion of the disc. We analyze existence and uniqueness for this problems in terms of the combinatorial structure of a word assigned to $f$. Our techniques are based on ideas of Blank who studied the extension problem in case of a planar target.

\end{abstract}

\section{Introduction}
\label{intro}
The Riemann Mapping Theorem and the Plateau problem are examples of extension problems.
The present work studies the following extension problem for immersions:
Suppose $f\colon \S^1 \to \S^2$ is a \emph{normal immersion}, that is, all self-intersection points are transversal double points. Does it extend to an immersion $F\colon\overline D \to \S^2$ from the closed unit disc to the sphere such that $F|_{\S^1}=f$?
A figure-$8$ curve shows that we cannot expect existence in general; the Milnor curve \cite[p. 11]{Poenaru} indicates that uniqueness (up to diffeomorphism) cannot be expected.\\
\\
The present work improves and extends an approach first introduced in 1967 by Samuel J. Blank in his PhD-thesis \cite{Blank}. His idea was to glue several embeddings together to obtain the desired immersion. Each embedding extends to a disc and hence these extensions can be glued together to an extension of the normal immersion. For that purpose a word is assigned to a normal immersion. By analyzing the structure of the word Blank was able to decide whether such a decomposition into embeddings could be achievied or not.

He showed that if the word contains only two special structures, called subwords, the corresponding normal immersion admits an extension to an immersed disc. Furthermore, different choices of the subwords result in different extensions.\\ 
\\
To extend this approach to the case of normal immersions into the sphere we introduce a new class of subwords (Section \ref{sec:groupings}). Moreover we represent different extensions in terms of wheigted trees. The tree structure shows that some of the steps of Blank are in fact not necessary. Therefore the present work not only extends Blank's results to a full classification of extensions of normal immersions  $f\colon\S^1 \to \S^2$ to immersed discs but it also improves the understanding by simplifying the proof. \\
\\
The PhD-thesis by Blank \cite{Blank} has never been published, perhaps since it is rather sketchy. Some of the missing proofs for the extension problem into the plane were later filled in by Valentin Po\'enaru \cite{Poenaru}. 

In 1974 George K. Francis \cite{Francis2} stated a classification result for surfaces immersed into the sphere. It provides a classification of immersed discs in the sphere, but the general results seems not to be appeared. Francis employs Blank's approach but develops it in a direction somewhat different from the present work. Let us also note that the results of the present work generalize to a complete classification of immersed surfaces in arbitrary closed surfaces \cite{Frisch}.\\
\\
We mention some more general problems which would be worth studying. Raising dimensions one might ask whether a normal immersion $f\colon\S^2 \to \S^3$ admits an extension to an immersed $3$-ball. 
A necessary condition for $f$ to extend is that the intersection with any preimage of a plane extends to an immersed surface (possibly disconnected). The longterm goal would be the classification of immersed $n$-manifolds in closed $n$-manifolds.\\
\\
The present work is based on the PhD-thesis of the author \cite{Frisch}. It was motivated by a classification of constant mean curvature surfaces in terms of discs immersed in $\S^2$, given by Grosse-Brauckmann, Kusner and Sullivan \cite{GKS}. So our results can also interpreted as existence and uniqueness results for constant mean curvature surfaces with finite topology.

\section{Basic  Topological Facts}
\label{sec:basics}

The present work studies the existence and uniqueness of extensions of normal immersions $f\colon\S^1 \to \S^2$ to immersed disc. The following results are true for more general immersions as well. Therefore the results are formulated for normal immersions $f\colon \coprod_{j=1}^m \S^1 \to N$ from the disjoint union of circles to arbitrary closed surface.

\begin{figure}
\psfrag{X1}{$X_1$} \psfrag{X2}{$X_2$} \psfrag{X3}{$X_3$} \psfrag{X4}{$X_4$}
\psfrag{k1}{\footnotesize $\psi(X_1)=1$} \psfrag{k2}{\footnotesize $\psi(X_2)=2$} 
\psfrag{k3}{\footnotesize $\psi(X_3)=0$} \psfrag{k4}{\footnotesize $\psi(X_4)=0$}
\psfrag{c}{$s$}
\includegraphics[width=\textwidth]{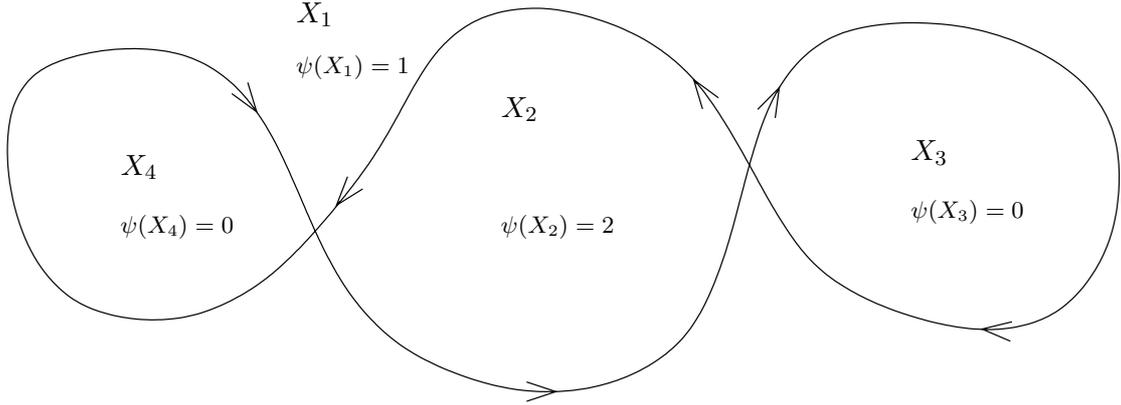}
\caption{Example of a normal numbering.}
\label{fig:numbering_ex}
\end{figure}

\begin{defs}
\emph{Suppose $N$ is a closed, orientable surface and $f\colon \coprod_{j=1}^m \S^1 \to N$ a normal immersion. Denote the components of $N \backslash f(\coprod_{j=1}^m \S^1)$ by $X_1,\ldots,X_k$. A normal immersion $f$ is called \emph{numerable} if a continuous function $\psi\colon \{X_1,\ldots,X_k\} \to \Z$ exists, such that 
\[
\psi(X_j)=\psi(X_i)+1,
\]
if $X_j$ and $X_i$ shares a common boundary and $X_j$ lies on the left of this boundary. \\
The function $\psi$ is called a \emph{numbering} of $f$. A function $\psi_n\colon \{X_1,\ldots,X_k\} \to \N_0$ with \linebreak $0 \in \im(\psi_n)$ is called a \emph{normal numbering} of $f$ (see Figure \ref{fig:numbering_ex}).}
\label{def:numerable}
\end{defs}

\begin{figure}
\centering
\psfrag{Xi}{$X_i$} \psfrag{Xj}{$X_j$} \psfrag{Xn}{$X_n$} \psfrag{Xm}{$X_m$}
\psfrag{ki}{\footnotesize $k_i$} \psfrag{kj}{\footnotesize $k_j=k_i-1$} 
\psfrag{kn}{\footnotesize $k_n=k_i-1$} \psfrag{km}{\footnotesize $k_m=k_i-2$}
\subfigure[]{\includegraphics[width=5cm]{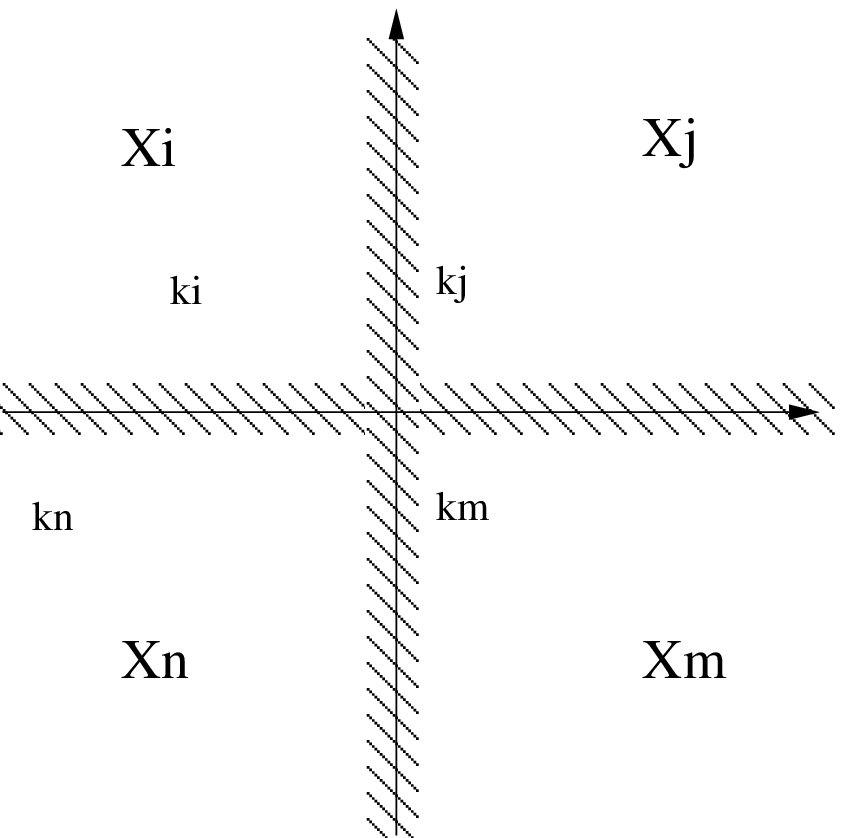}}
\qquad
\psfrag{Xi}{$X_i$} \psfrag{Xj}{$X_j$} \psfrag{Xn}{$X_n$} \psfrag{Xm}{$X_m$}
\psfrag{ki}{\footnotesize $k_i$} \psfrag{kj}{\footnotesize $k_j=k_i+1$} 
\psfrag{kn}{\footnotesize $k_n=k_i-1$} \psfrag{km}{\footnotesize $k_m=k_i$}
\subfigure[]{\includegraphics[width=5cm]{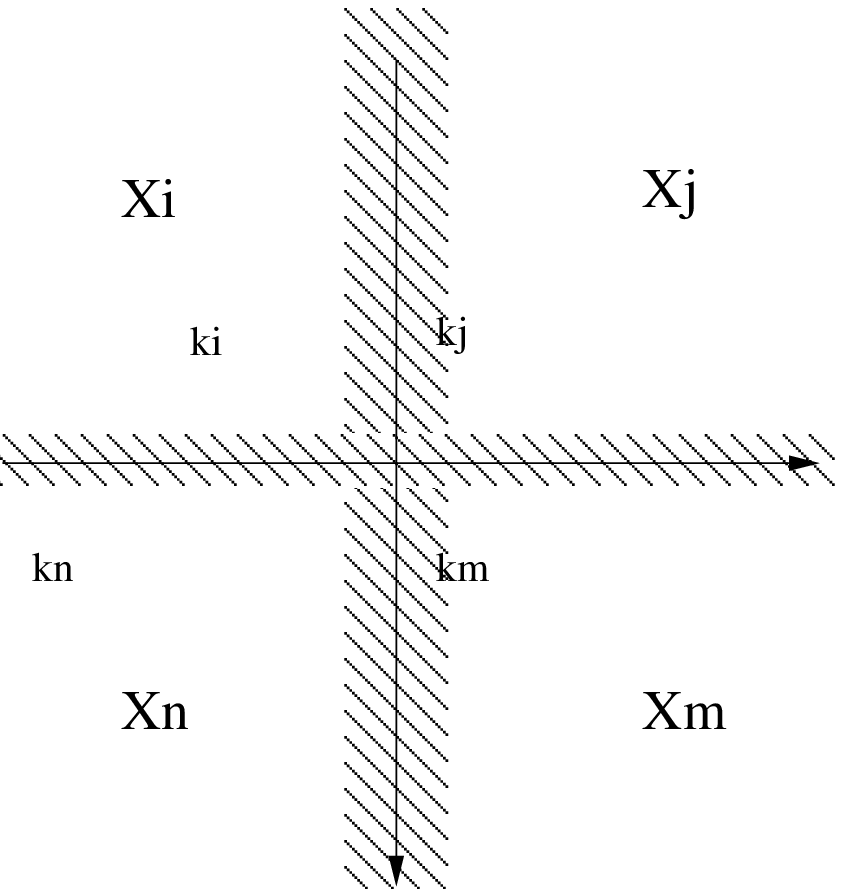}}
\caption{The path shown travels from left to right and crosses a double point. If the two components on the left are numbered as indicated, the numbering extends to the components on the right as indicated.}
\label{fig:algo_intersect}
\end{figure}

In \cite{MC} Margaret McIntyre and Grant Cairns describe an algorithm that assigns a numbering to a given normal immersion $f\colon \S^1 \to \R^2$. That algorithm can be applied to general normal immersions $f\colon \coprod_{j=1}^m \S^1 \to N$ as well and is shown in Figure \ref{fig:algo_intersect}.
It  assigns a numbering to each numerable immersion \linebreak $f\colon\coprod_{j=1}^m\S^1 \to N$. From this numbering we can derive a normal numbering by setting
\begin{equation}
\psi_n(X_l):=\psi(X_l)-\min_{i=1,\ldots,k} \psi(X_i).
\label{eq:normal_numbering}
\end{equation}

We will now analyze when a normal immersion is numerable. For that suppose $N$ is a closed surface of genus $g$ and let $\nu_1,\ldots,\nu_{2g}$ be a set of differentiable curves whose homology classes $[\nu_1],\ldots,[\nu_{2g}]$ form a generating set for the first homology group $H_1(N)$. 
Now let $f$ be an oriented, normal curve with homology class
\[
[f]=n_1[\nu_1]+\ldots+n_{2g}[\nu_{2g}].
\]
Moreover assume that the curves $f,\nu_1,\ldots,\nu_{2g}$ are pairwise transversal.

\begin{lemma}[McIntyre, Cairns 1993]
Let $X_1,\ldots,X_k$ be the components of \linebreak$N\backslash\{f \cup \nu_1 \cup \ldots \nu_{2g}\}$. One can associate integers to each of the components $X_1,\ldots,X_k$ such that at each segment of $f$ the number to the left of $f$ is $1$ greater than the number to the right of $f$, and for each $i=1,\ldots,2g$ the number to the left of each segment of $\nu_i$ is $n_i$ less than the number to the right of $\nu_i$.

The numbering is unique if we choose one of the components to be $0$.
\label{lem:MC}
\end{lemma}

\begin{proof}
See \cite{MC}, Lemma 2.
\end{proof}

This yields the following lemma:

\begin{lemma}
Let $M$ be a compact connected oriented manifold with boundary \linebreak $\partial M=\coprod_{j=1}^m \S^1$. Let $N$ be a closed surface and $f\colon \coprod_{j=1}^m \S^1 \to N$ a normal immersion. Then $f$ is numerable if and only if $f$ is homologically trivial. 
\label{lem:homolog}
\end{lemma}

\begin{proof}
See \cite{Frisch}, Lemma 1.2.5.
\end{proof}

\begin{cor}
If $N$ is simply connected then each normal immersion is numerable.
\label{cor:numerable}
\end{cor}

If $f\colon\coprod_{j=1}^m\S^1 \to N$ extends to an immersion $F\colon M \to N$, the value of a normal numbering  marks the difference between the number of preimages $\omega_f(X_j)$ of a component $X_j$ under $f$.
Indeed each time we pass the curve $f(\coprod_{j=1}^m \S^1)$ from left to right, we add a layer and each time we pass it from right to left we lose a layer of the surface (recall that the surface lies to the left of the curve $f(\coprod_{j=1}^m \S^1)$). In Figure \ref{fig:algo_intersect} the surface is marked by the pattern. \\
\\
The number of preimages of an extension $F\colon M \to N$ defines a numbering of the boundary curve $f\colon \coprod_{j=1}^m \S^1 \to N$. While the number of preimages of different extensions is not constant in general, it is for $N=\S^2$ :

\begin{prop}
Let $f\colon \S^1 \to \S^2$ be a normal immersion which extends to $\overline D$ and denote the components of $\S^2 \backslash f(\S^1)$ by $X_1,\ldots,X_k$. The number of preimages of an extension $F\colon\overline D \to \S^2$ for each $X_j$ is uniquely determined by $f$.
\label{prop:preimage_const}
\end{prop}

\begin{proof}
The proof is based on the application of the global Gauss-Bonnet Theorem. See \cite{Frisch}, Proposition 1.3.1. 

\end{proof}

That satisfies the following definition:

\begin{defs}
\emph{Let $f\colon\S^1\to\S^2$ be a normal immersion and $\psi_n$ a normal numbering of $f$. For each component $X_j$ of $\S^2 \backslash f(\S^1)$ the number}
\[
\omega_j:=\psi_n(X_j)
\]
\emph{is called the \emph{degree of $X_j$}\index{degree!of a component}.}
\end{defs}

\begin{figure}
\centering
\psfrag{f1}{$f(t_1)=f(t_2)$}
\psfrag{f2}{$f'(t_1)=f'(t_2)$}
\psfrag{f3}{$f_1$}
\psfrag{f4}{$f_2$}
\subfigure[]{\includegraphics[width=.3\textwidth]{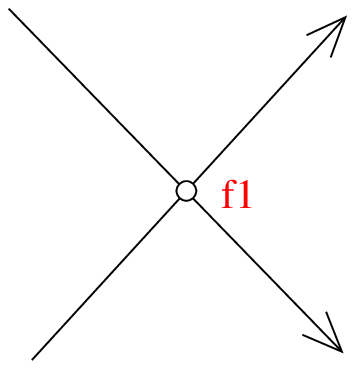}}\quad\subfigure[]{\includegraphics[width=.3\textwidth]{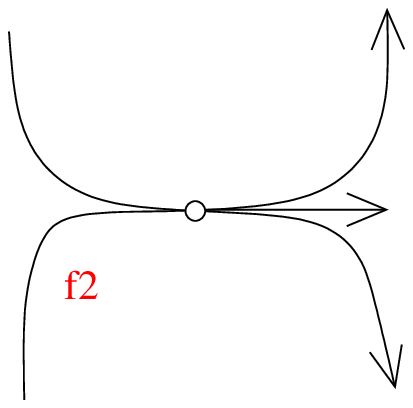}}\quad\subfigure[]{\includegraphics[width=.3\textwidth]{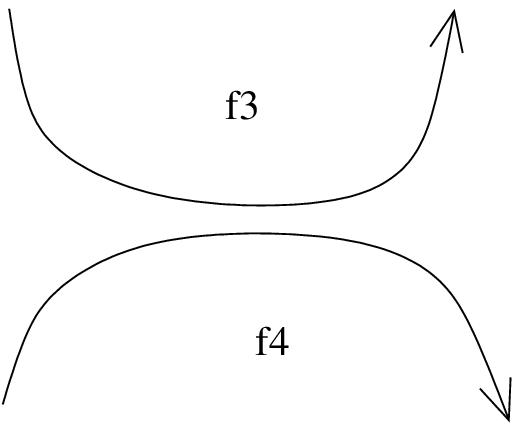}}
\caption{Decomposing a normal immersion into embeddings}
\label{fig:embedding}
\end{figure}

As a conclusion of the proof we get a necessary condition for a normal immersion to extend in terms of the tangent winding number. The definition of the tangent winding number is based on a decomposition of the normal immersion into embeddings. To decompose the normal immersion adjust the tangent vectors in the double points and afterwards remove the double points (see Figure \ref{fig:embedding}).

\begin{defs}
\emph{Let $f\colon \S^1 \to \S^2$ be a normal immersion and $f_1,\ldots,f_n$ the family of embeddings which results from the decomposition. Choose a base point $x_0 \in \S^2 \backslash f(\S^1)$ in a component of minimal degree. With an orientation preserving diffeomorphism $\phi\colon\S^2 \backslash \{x_0\} \to \R^2$ we denote by $\widetilde f_i=\phi \circ f_i$ the induced embedding to $\R^2$.}

\emph{Define the \emph{tangent winding number} of $f_i$ as }
\[
\tau(f_i):=\wind(\widetilde f_i)
\]
\emph{and the \emph{tangent winding number of $f$} as}
\[
\tau(f)=\sum_{i=1}^n \tau(f_i).
\]
\label{def:rank}
\end{defs}

\begin{prop}[Necessary Extension Condition for one Boundary Component]
\index{Proposition!Necessary Extension Condition!for one boundary component}
Every normal immersion $f\colon \S^1 \to \S^2$ which can be extended to an immersion $F\colon \overline D \to \S^2$ has the tangent winding number
\[
\tau(f)=1-2\omega_1 \qquad \textrm{for } \omega_1 \in \N_0,
\]
where $\omega_1$ denotes the degree of the base point component $X_1$.
\label{prop:notwendig}
\end{prop}

\section{Groupings}
\label{sec:groupings}

In this section we describe how to use Blank's algorithm \cite{Blank} in the case of an immersion $f\colon \S^1 \to \S^2$ to define a word $w(f)$. The word will depend on some choices. A choice optimal for our purposes will lead to what we call a reduced word.\\
\\
To get the word choose a component $X_j$ with minimal degree of $\S^2 \backslash f(\S^1)$ as base component and pick a base point $x_0$ in $X_j$. Pick for each other component $X_i$ with $i \neq j$ a point $p_i$. Then each $p_i$ can be connected with the base point $x_0 \in X_j$ by a ray $\hat c_i$ such that
\begin{itemize}
\item $\hat c_i \cap \hat c_l = \{x_0\}$ for $i \neq l$ and
\item each $\hat c_i$ has a minimal number of intersection points with $f(\S^1)$.
\end{itemize}

\begin{figure}
\centering
\psfrag{x0}{{$x_0$}}
\psfrag{s}{{$s$}}
\psfrag{p1}{{$p_1$}}
\psfrag{p2}{{$p_2$}}
\psfrag{p3}{{$p_3$}}
\psfrag{p4}{{$p_4$}}
\psfrag{p5}{{$p_5$}}
\psfrag{c1}{{$\hat c_1$}}
\psfrag{c2}{{$\hat c_2$}}
\psfrag{c3}{{$\hat c_3$}}
\psfrag{c4}{{$\hat c_4$}}
\psfrag{c5}{{$\hat c_5$}}
\psfrag{a1+}{{$a_1$}}
\psfrag{a2+}{{$a_2$}}
\psfrag{a3+}{{$a_3$}}
\psfrag{a4+}{{$a_4$}}
\psfrag{a5+}{{$a_5$}}
\psfrag{a1-}{{$a_1^{-1}$}}
\psfrag{a2-}{{$a_2^{-1}$}}
\psfrag{a3-}{{$a_3^{-1}$}}
\psfrag{a4-}{{$a_4^{-1}$}}
\psfrag{a5-}{{$a_5^{-1}$}}
\includegraphics[width=.49\textwidth]{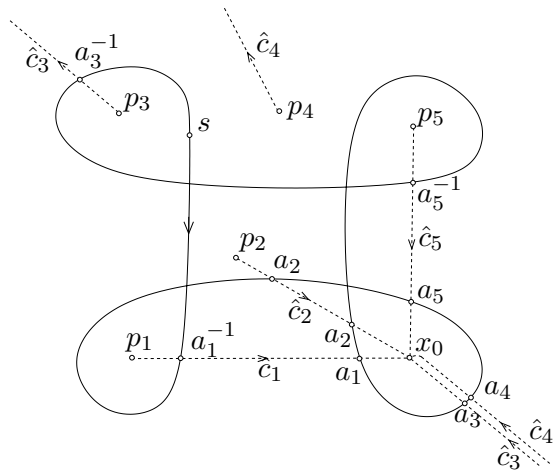}
\caption{The system of rays and the intersection points in projection on $\R^2$. The word in that case is $a_1^{-1}a_2a_5a_4a_3a_1a_2a_5^{-1}a_3^{-1}$.}
\label{fig:ex1}
\end{figure}

Label an intersection point of the ray $\hat c_i$ and $f(\S^1)$ with $a_i$ of $\hat c_i$ crosses from left to right and with $a_i^{-1}$ if it crosses from right to left. After labeling all intersection points write them down in the order they appear when you walk along $f(\S^1)$, starting at an initial point $s$ (see Figure \ref{fig:ex1}).

\begin{defs}
\emph{Let $f\colon \S^1 \to \S^2$ be a normal immersion. The ordered collection of intersection points is called \emph{word of $f$}\index{word!of $f$} and will be denoted by $w(f)$. It is unique up to cyclic permutations.}
\label{def:word}
\end{defs}

\begin{bem}
\emph{The uniqueness up to cyclic permutations is caused by the choice of the initial point $s$. The choice of a different starting point leads to a cyclic permutation of the word, since the order of the intersection points does not change. Only the choice of the first letter is different.}
\end{bem}

\begin{defs}
\emph{A \emph{subword}\index{word!subword}\index{subword|see{word}} $\omega$ of $w(f)$ is a subsequence of consecutive letters of $w(f)$.} 
\emph{A word $w(f)$ is called} reduced\index{word!reduced} \emph{if it satisfies the following conditions:}

\begin{itemize}
\item
\emph{The word $w(f)$ contains no subword of the form $a_ia_i^{-1}$ or $a_i^{-1}a_i$.}

\item
\emph{The word $w(f)$ is not of the form $a_i \ldots a_i^{-1}$ or $a_i^{-1}\ldots a_i$.}
\end{itemize}
\label{def:reduce}
\end{defs}

The reduced word is unique up to cyclic permutations as well. If we move the starting point $s$ along $f(\S^1)$ the reduced word changes by cyclic permutations. 

We will now show that for each normal immersion $f\colon \S^1 \to \S^2$ a word $w(f)$ exists which is reduced:

\begin{lemma}
For every normal immersion $f\colon \S^1 \to \S^2$ the points $p_i$, $s$ and $x_0$ as well as the rays $\hat c_i$ can be chosen such that the resulting word $w(f)$ is reduced.
\end{lemma}

\begin{proof}
Assume the word $w(f)$ of $f$ is not reduced.

\begin{enumerate}
\item Assume $w(f)$ contains a subword of the form $a_ia_i^{-1}$. That means that $f(\S^1)$ hits the ray $\hat c_i$ two times consecutively. Since there is no other letter between $a_i$ and $a_i^{-1}$ the situation looks like in Figure \ref{fig:reduce1} (a). In this case we can homotopy the ray $\hat c_i$ such that the intersection points $a_i$ and $a_i^{-1}$ do not occur any more (Figure \ref{fig:reduce1} (b)).
If $w(f)$ contains a subword of the form $a_i^{-1}a_i$ then the situation is the same. Only the orientation of the ray $\hat c_i$ in relation to $f(\S^1)$ has changed. But there is still no intersection of another ray and $f(\S^1)$ between $a_i^{-1}$ and $a_i$ and so we can homotopy $\hat c_i$ again such that the letters $a_i^{-1}a_i$ do not occur any more.

\begin{figure}
\centering
\psfrag{f}{$f(\S^1)$}
\psfrag{ci}{{$\hat c_i$}}
\psfrag{ai+}{{$a_i$}}
\psfrag{ai-}{{$a_i^{-1}$}}
\subfigure[Situation when $w(f)$ is not reduced.]{\includegraphics[width=.45\textwidth]{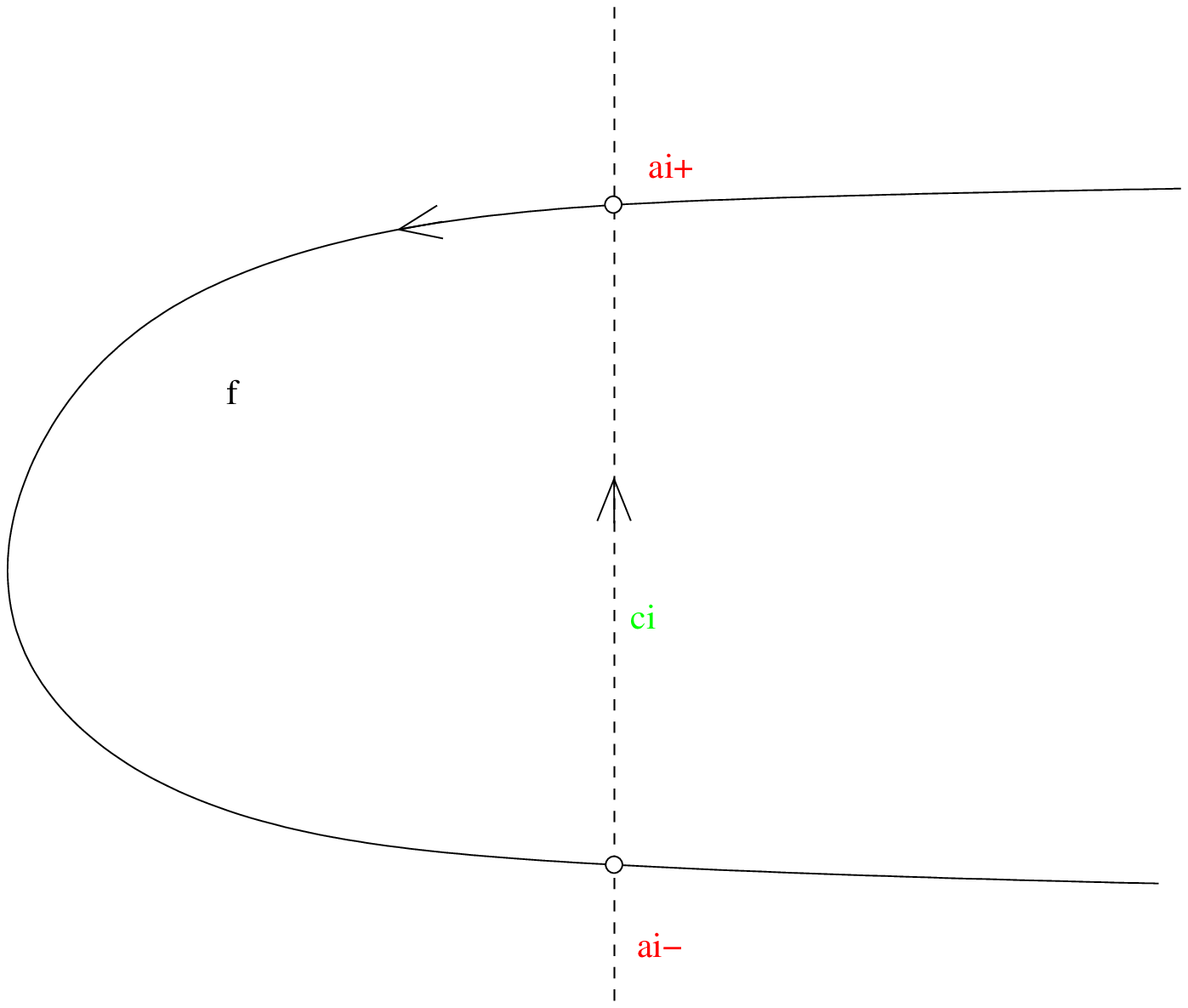}}\qquad
\subfigure[Changing $\hat c_i$ such that the resulting word is reduced.]{\includegraphics[width=.45\textwidth]{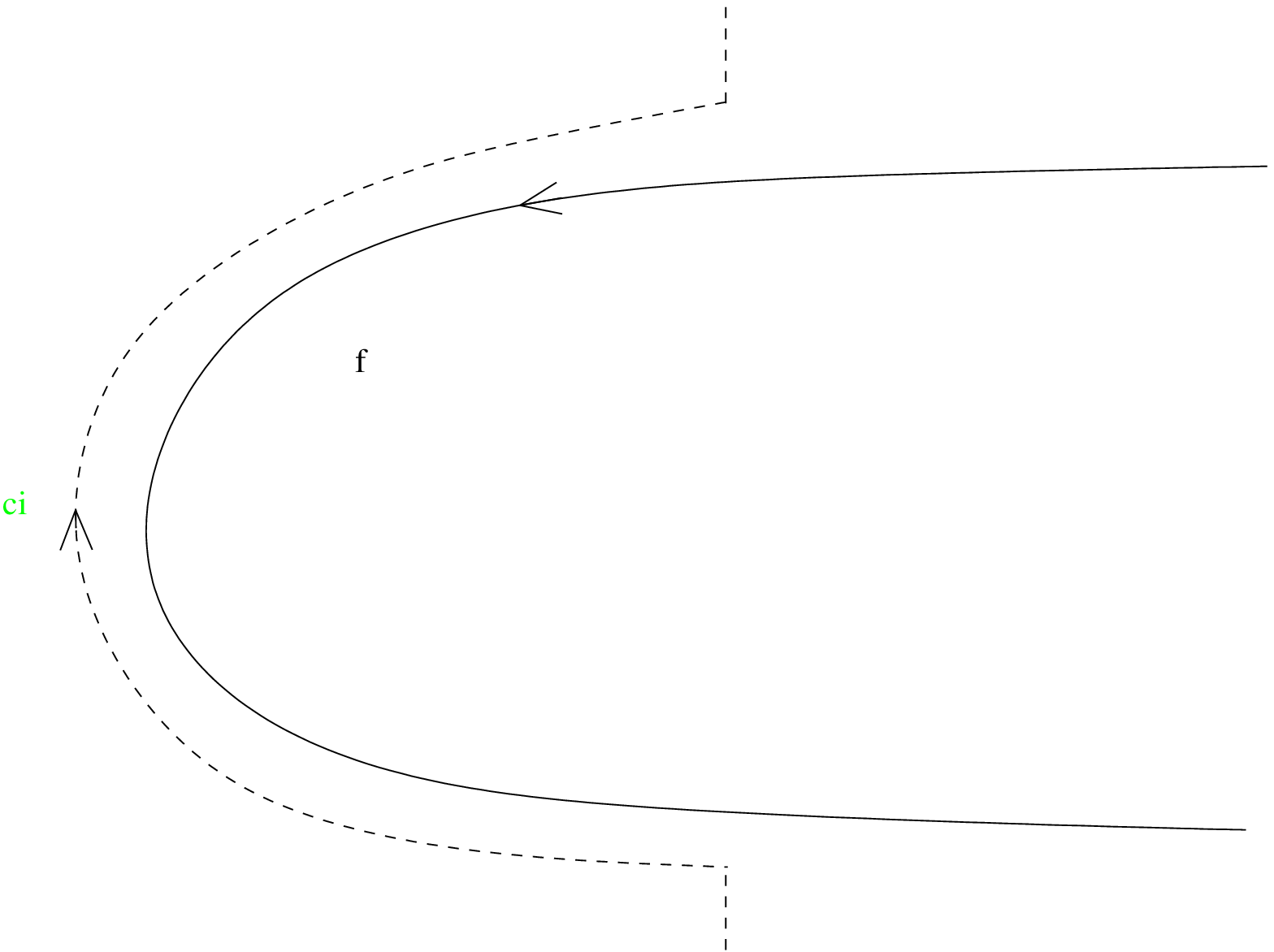}}
\caption{Reducing the word $w(f)$.}
\label{fig:reduce1}
\end{figure}

\item Assume the word $w(f)$ is of the form $a_i^{-1} \ldots a_i=a_i^{-1}\omega a_i$, where $\omega$ is the subword which contains all letters except the outer ones. That means, that the initial point $s$ lies between $a_i$ and $a_i^{-1}$. Choose a new initial point $\widetilde s$ which is directly after the intersection point $a_i^{-1}$. Now we get a word $\omega a_ia_i^{-1}$. As in $1.$ this word can be reduced to $\omega$. 
\end{enumerate}

This shows that every time one of the conditions of Definition \ref{def:reduce} is not satisfied we can change the rays or the initial point to get a setting for which the word is reduced.
\end{proof}

In the following the word $w(f)$ of a normal immersion $f$ is assumed to be reduced. As seen in  Proposition \ref{prop:notwendig} a necessary condition for a normal immersion $f\colon \S^1 \to \S^2$ to extend to $\overline D$ is $\tau(f)=1-2\omega_1$. 

A sufficient condition for $f$ to extend is provided by the combinatorial structure of the word. Crucial for that condition are special subwords which we will introduce first. 

\begin{defs}
\emph{For a reduced word $w(f)$ the subword}
\begin{itemize}
\item
\emph{$a_ia_j\ldots a_k$ is called a \emph{positive word}\index{word!positive word}\index{positive word|see{word}} if consecutive letters are different,}

\item
\emph{$a_i^{\pm 1} p a_i^{\mp 1}$ is called a \emph{pairing}\index{word!pairing}\index{pairing|see{word}} if $p$ is a positive word,}

\item
\emph{$a_j^{-1}a_i^{-1}$ is called a \emph{negative group}\index{word!negative group}\index{negative group|see{word}} if $i\neq j$.}
\end{itemize}

\emph{The empty word is defined as positive.}
\label{def:subword}
\end{defs}

We will show that an immersion whose reduced word is composed of these special subwords extends. If a reduced word $w(f)$ contains a subword $\omega$, which is a pairing or a negative group, $w(f)$ can be written as $w(f)=x_1\omega x_2$. Then we can cancel the subword $\omega$ out of $w(f)$ to get $\widetilde{w}(f)=x_1x_2$. 

\begin{defs}
\emph{A reduced word $w(f)$ is called \emph{groupable}\index{word!groupable}\index{groupable|see{word}} if a cancellation of pairings and negative groups exists, such that a positive word remains. }
\end{defs}

\begin{example}
\emph{A groupable word which can be cancelled in two different ways:}
\begin{eqnarray*}
a_2{\bf |a_3^{-1}a_1^{-1}|}a_4^{-1}a_2a_1a_4a_3 &\leadsto a_2{\bf |a_4^{-1}a_2a_1a_4|}a_3 &\leadsto a_2a_3 \\
a_2a_3^{-1}{\bf |a_1^{-1}a_4^{-1}|}a_2a_1a_4a_3 &\leadsto a_2{\bf |a_3^{-1}a_2a_1a_4a_3|} &\leadsto a_2
\end{eqnarray*}
\label{bsp:grouping}
\end{example}

\begin{figure}
\centering
\psfrag{x0}{{$x_0$}}
\psfrag{ai+}{{$a_i$}}
\psfrag{ai-}{{$a_i^{-1}$}}
\psfrag{aj-}{{$a_j^{-1}$}}
\psfrag{ak-}{{$a_k^{-1}$}}
\psfrag{ci}{{$\hat c_i$}}
\psfrag{cj}{{$\hat c_j$}}
\psfrag{ck}{{$\hat c_k$}}
\psfrag{pi}{{$p_i$}}
\psfrag{pj}{{$p_j$}}
\psfrag{pk}{{$p_k$}}
\psfrag{[]}{{$[a_i^{-1},a_i]$}}
\psfrag{()}{{$[a_j^{-1},a_k^{-1}]$}}
\includegraphics[width=\textwidth]{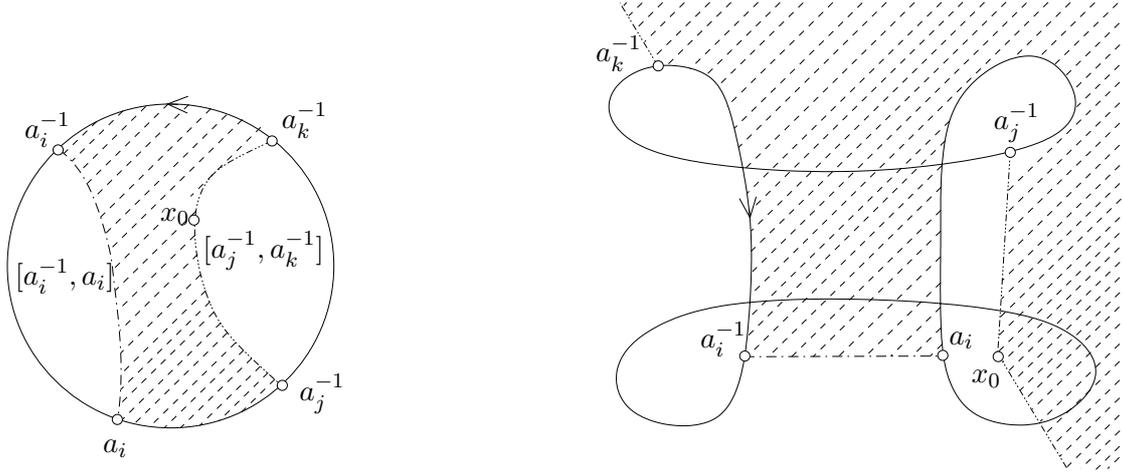}
\caption{Extending the normal immersion to intervals. The pattern marks the disc enclosed by the intervals.}
\label{fig:interval}
\end{figure}

The goal is to show that a normal immersion $f$ whose reduced word $w(f)$ is groupable extends to an immersion on $\overline D$ and vice versa. \\
\\
Suppose $f\colon \S^1 \to \S^2$ is a normal immersion and denote the corresponding reduced word by $w(f)$. We mark the preimages of the intersection points $a_i^{\pm 1}$ on $\S^1$. For convenience we use the notation $a_i^{\pm 1}$ for the preimages as well.

Given the case of a pairing $a_i^{\pm 1}\ldots a_i^{\mp 1}$ we can join the corresponding points in the preimage by a differentiable curve in $D$ which is transverse to $\S^1$. This curve is called an interval and is denoted by $[a_i^{\pm 1},a_i^{\mp 1}]$ (Figure \ref{fig:interval}).\\
\\
In the case of a negative group $a_j^{-1}a_k^{-1}$, we join the points in the preimage by an interval which contains a preimage of the base point $x_0$. Of course, this can be the case only if a possible extension $F\colon \overline D\to \S^2$ is surjective. That means the base point component has a positive degree. In that case mark an arbitrary point in the interior of $\overline D$ as $x_0$. Now join $a_j^{-1}$ and $x_0$ by an interval as well as $x_0$ and $a_k^{-1}$. The intervals should be chosen as transversal to $\S^1$ and differentiable in $x_0$. Denote the union of these two intervals by $[a_j^{-1},a_k^{-1}]$ (Figure \ref{fig:interval}).\\
\\

If there is more than one such interval the question arises whether these intervals can be chosen, so that they are disjoint.

\begin{lemma}
Suppose $f\colon \S^1 \to \S^2$ is a normal immersion with groupable word $w(f)$. Then for each pairing and each negative group an interval can be chosen such that all intervals are disjoint. 
\label{lem:disjoint}
\end{lemma}

\begin{proof}
See \cite{Frisch}, Lemma 2.4.5.
\end{proof}

\begin{figure}
\centering
\psfrag{I1}{{$I_1$}}
\psfrag{I2}{{$I_2$}}
\psfrag{I3}{{$I_3$}}
\psfrag{I4}{{$I_4$}}
\psfrag{a1}{{$\overset{\textrm{No. }1}{a_1^{-1}}$}}
\psfrag{b1}{{$a_2^{-1}$}}
\psfrag{a2}{{$a_3^{-1}$}}
\psfrag{b2}{{$a_4^{-1}$}}
\psfrag{a3}{{$a_1$}}
\psfrag{b3}{{$a_1^{-1}$}}
\psfrag{a4}{{$a_4^{-1}$}}
\psfrag{b4}{{$a_4$}}
\psfrag{a5}{{$a_1^{-1}$}}
\psfrag{b5}{{$a_3^{-1}$}}
\psfrag{D1}{{$D_1$}}
\psfrag{D2}{{$D_2$}}
\psfrag{D3}{{$D_3$}}
\psfrag{D4}{{$D_4$}}
\psfrag{D5}{{$D_5$}}
\psfrag{1}{{$1$}}
\psfrag{2}{{$2$}}
\psfrag{3}{{$3$}}
\psfrag{4}{{$4$}}
\psfrag{5}{{$5$}}
\psfrag{6}{{$6$}}
\psfrag{7}{{$7$}}
\psfrag{8}{{$8$}}
\psfrag{9}{{$9$}}
\psfrag{10}{{$10$}}
\psfrag{12}{{$12$}}
\psfrag{34}{{$34$}}
\psfrag{510}{{$510$}}
\psfrag{69}{{$69$}}
\psfrag{78}{{$78$}}
\subfigure[The tree]{\includegraphics[width=.4\textwidth]{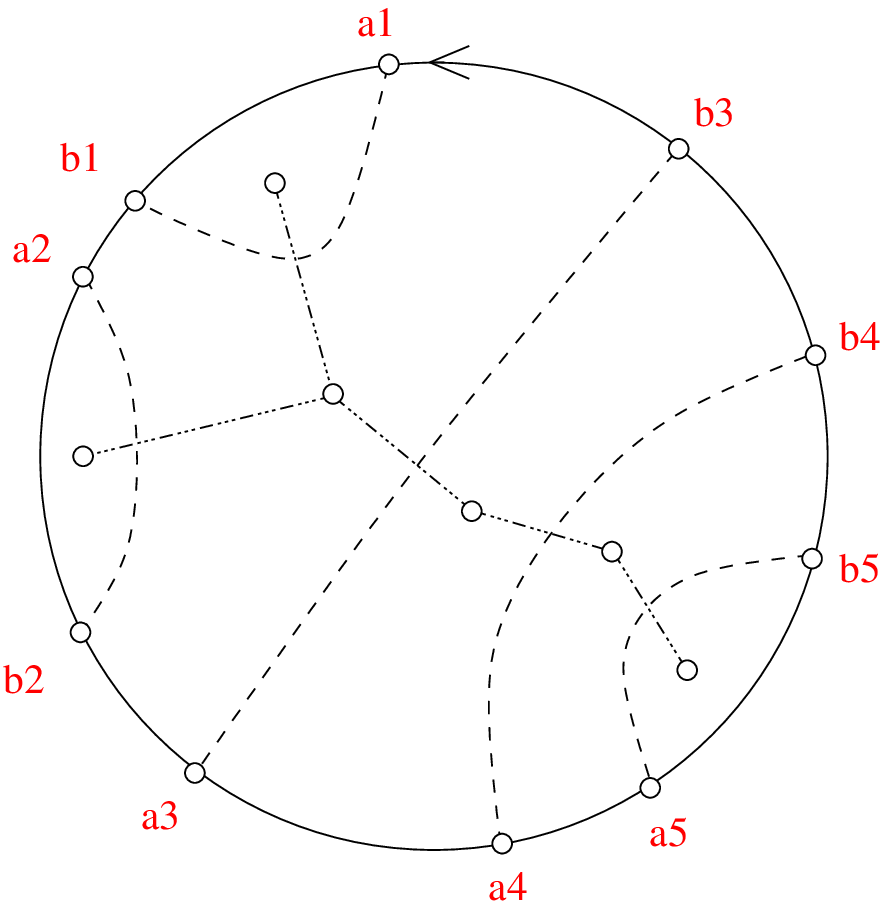}}\quad
\subfigure[The weighted tree]{\includegraphics[width=.4\textwidth]{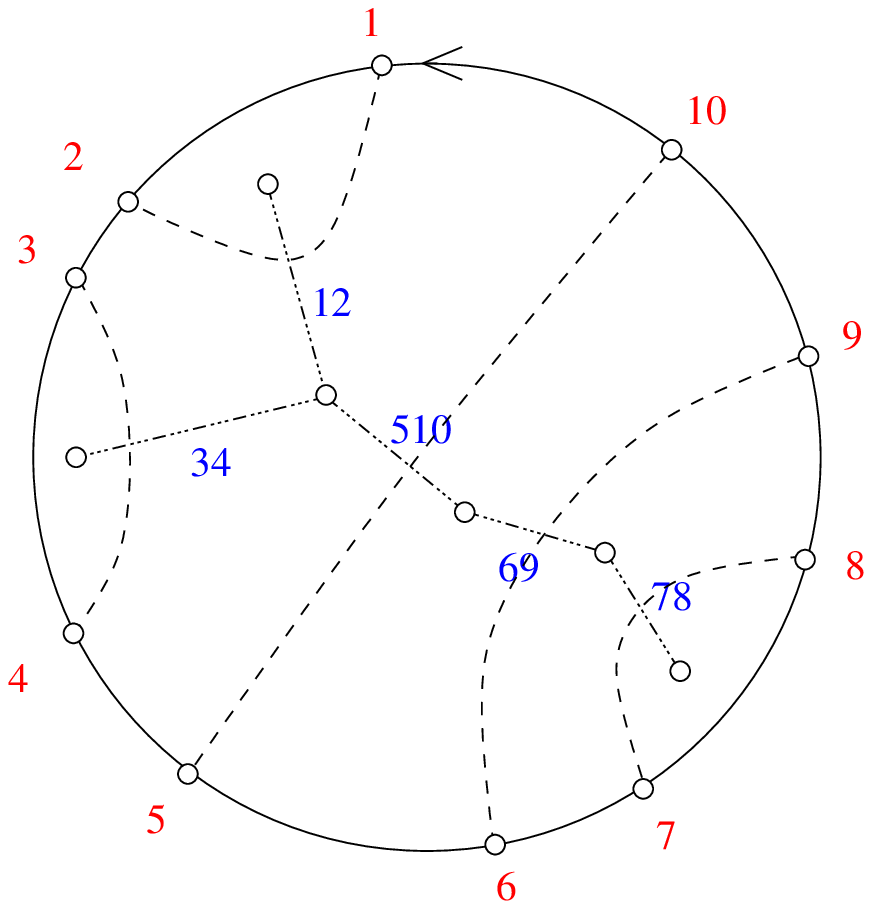}}
\caption{Decomposing $\overline D$ into smaller discs.}
\label{fig:small_disks}
\end{figure}

Choosing the intervals to be disjoint leads to a subdivision of $\overline D$ into smaller discs. If we put a vertex in each small disc and connect two vertices of adjacent discs, we get a tree which is dual to $\S^1$ together with the intervals.
Take the reduced word $w(f)$ and number all letters. Remember that the reduced word is unique up to cyclic permutations so the numbering is unique if we pick a first letter. 

Now mark the preimages of the letters and their number in $\S^1$. Since $w(f)$ is groupable there has to be a pairing or a negative group which can be canceled. Cancel it and insert the corresponding interval into $\overline D$. Continue until a positive word remains. Now $\overline D$ is decomposed into smaller discs (Figure \ref{fig:small_disks}). \\
\\
Put in each of these discs a vertex and connect two vertices if their small discs share a common boundary. Since the intervals can be chosen so that they are disjoint (Lemma \ref{lem:disjoint}) each edge intersects exactly one interval $[a,b]$. 
If the letter $a$ has the number $i$ and the letter $b$ has the number $j$ then label the edge with $ij$ (Figure \ref{fig:small_disks} (b)).\\
\\
If $w(f)$ contains no pairing or negative group but is groupable then $w(f)$ is a positive word. In that case the decomposition consists of only one disc and hence the induced tree has only one vertex and no edge.

\begin{lemma}
Suppose $f:\S^1 \to \S^2$ is a normal immersion and the corresponding reduced word is $w(f)$. If $w(f)$ is groupable each cancellation, which results in a positive word, induces a weighted tree.
\label{lem:tree}
\end{lemma}

\begin{proof}
See \cite{Frisch}, Lemma 2.4.6.
\end{proof}

\begin{defs}
\emph{Suppose $f\colon \S^1 \to \S^2$ is a normal immersion and $w(f)$ the corresponding reduced word. If $w(f)$ is groupable then the induced weighted tree $\G$ is called a \emph{grouping of $f$}.}

\emph{Two groupings of $w(f)$ are \emph{equivalent}\index{grouping!equivalent} if the weighted trees are isomorphic. }
\label{def:grouping_different}
\end{defs}

\begin{figure}
\centering
\psfrag{a1-a4-}{{$23$}}
\psfrag{a5-a5}{{$46$}}
\psfrag{a4-a5-}{{$34$}}
\psfrag{a1-a1}{{$25$}}
\psfrag{a1}{{$5$}}
\psfrag{a2}{{$1$}}
\psfrag{a3}{{$a_3$}}
\psfrag{a4}{{$7$}}
\psfrag{a5}{{$6$}}
\psfrag{a1-}{{$2$}}
\psfrag{a4-}{{$3$}}
\psfrag{a5-}{{$4$}}
\subfigure[]{\includegraphics[width=.4\textwidth]{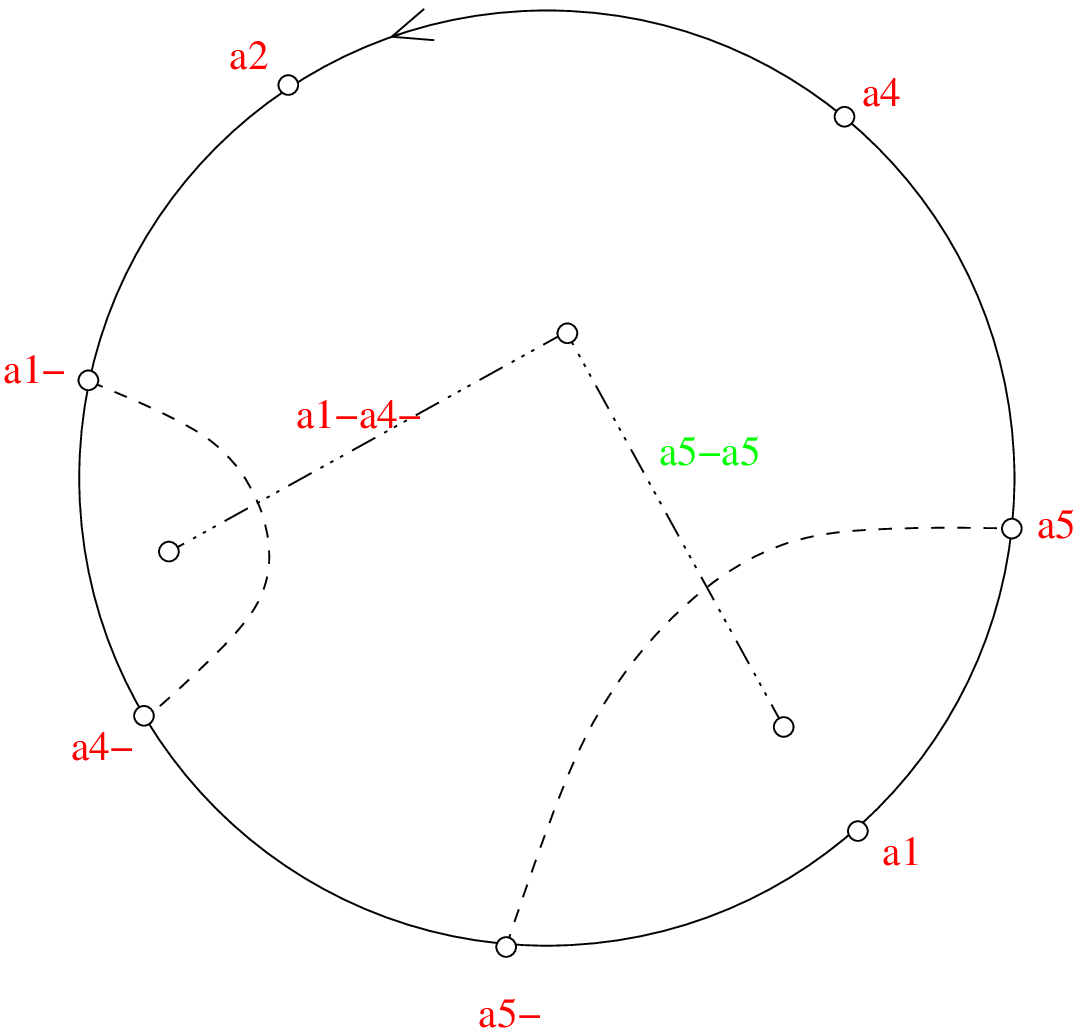}}\qquad
\subfigure[]{\includegraphics[width=.4\textwidth]{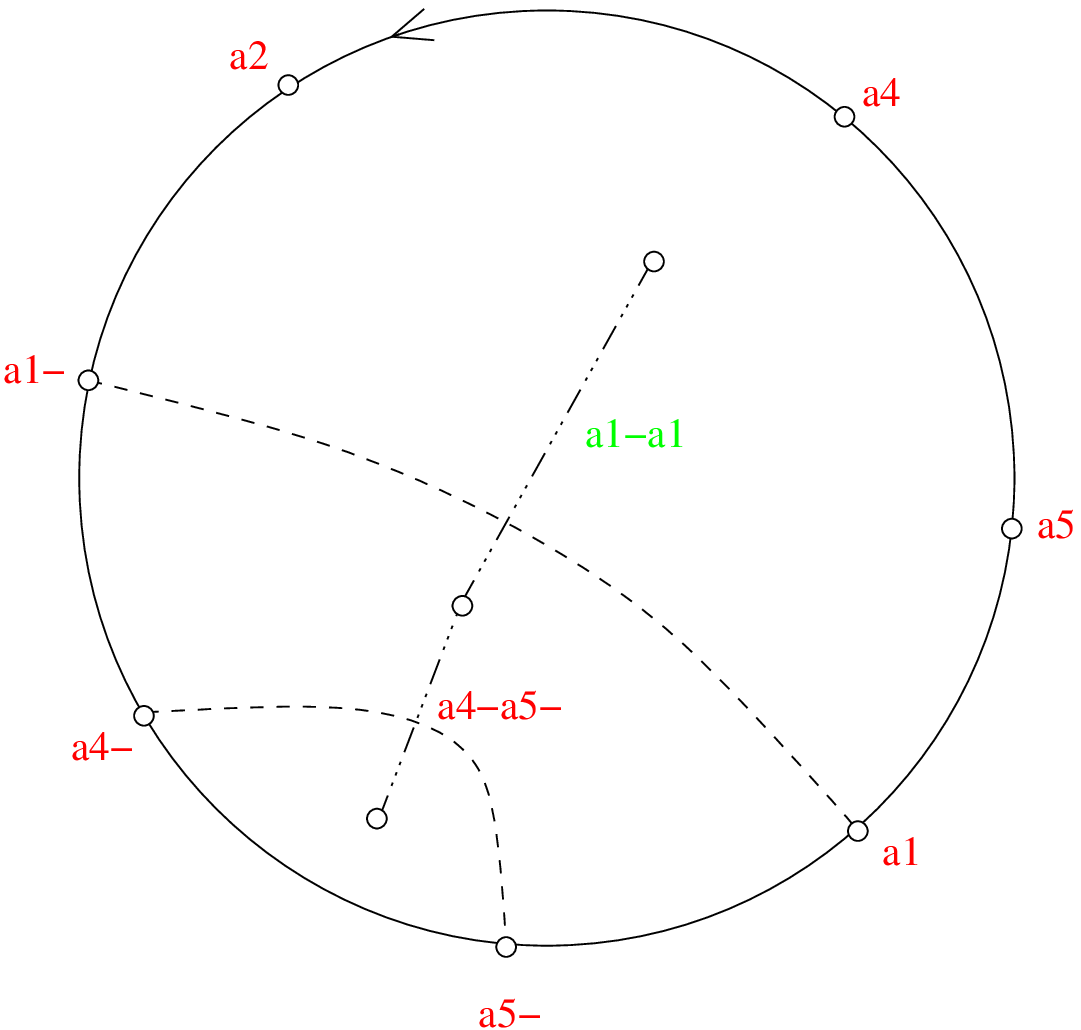}}
\caption{Two different groupings of the word $a_2a_1^{-1}a_4^{-1}a_5^{-1}a_1a_5a_4$}
\label{fig:cancel}
\end{figure}

\section{Immersed Discs into the Sphere}
\label{sec:disc->sphere}

After some preparations in the previous section we can now look at the main goal: Extending immersions.

\subsection{Existence of Extensions}

For convenience we use the following speech:

\begin{defs}
\emph{We say that a word $w(f)$} surrounds an immersed disc\index{word!surrounds an immersed disc} \emph{if there is a normal immersion $f\colon \S^1 \to \S^2$ which has the reduced word $w(f)$ and which can be extended to an immersion $F\colon \overline D \to \S^2$. 
If the same holds true for an embedding, $w(f)$ \emph{surrounds an embedded disc}\index{word!surrounds an embedded disc}.}

\emph{A subword $a\dots b$ of a word $w(f)$} surrounds an immersed disc\index{word!subword!surrounds an immersed disc} \emph{if there is an interval $[a,b]$, such that $\S^1_{[a,b]} \cup [a,b]$ surrounds an immersed disc.
Again, if the same holds true for an embedding then the subword \emph{surrounds an embedded disc}\index{word!subword!surrounds an embedded disc}.}
\end{defs}

We will show that a normal immersion $f$ extends to $\overline D$ if and only if the reduced word $w(f)$ is groupable. Since groupable means to cancel out pairings and negative groups in $w(f)$ we will show at first that these subwords themselves surround  immersed discs.

\begin{lemma}
A positive word surrounds an embedded disc.
\label{lem:positive_word}
\end{lemma}

\begin{proof}
Decompose the normal immersion $f\colon \S^1 \to \S^2$ into embeddings $f_1,\ldots,f_n$. W.l.o.g. $\tau(f_1)=+1$. Assume there is an embedding with $\tau(f_j)=-1$. Since $x_0$ lies in a component of minimal degree the ray $\hat c_j$ has to intersect $f_j$ from right to left. This produces a negative intersection point $a_j^{-1}$ in contradiction to the assumption that $w(f)$ is a positive word. Hence all embeddings have tangent winding number $\tau(f_j)=+1$ which leads to
\[
\tau(f)=\sum_{j=1}^n \tau(f_j)=n.
\]
On the other hand $\tau(f) \leq 1$, according to Proposition \ref{prop:notwendig}. Thus $n=1$ and hence $f$ is an embedding itself. Therefore it extends to an embedding $F\colon \overline D \to \S^2$ according to the Riemann Mapping Theorem (see \cite{Ahlfors}, p. 172ff).
\end{proof}

\begin{figure}
\centering
\psfrag{ci}{{$\hat c_i$}}
\psfrag{+}{{$+$}}
\psfrag{ai+}{{$a_i$}}
\psfrag{ai-}{{$a_i^{-1}$}}
\includegraphics[width=.9\textwidth]{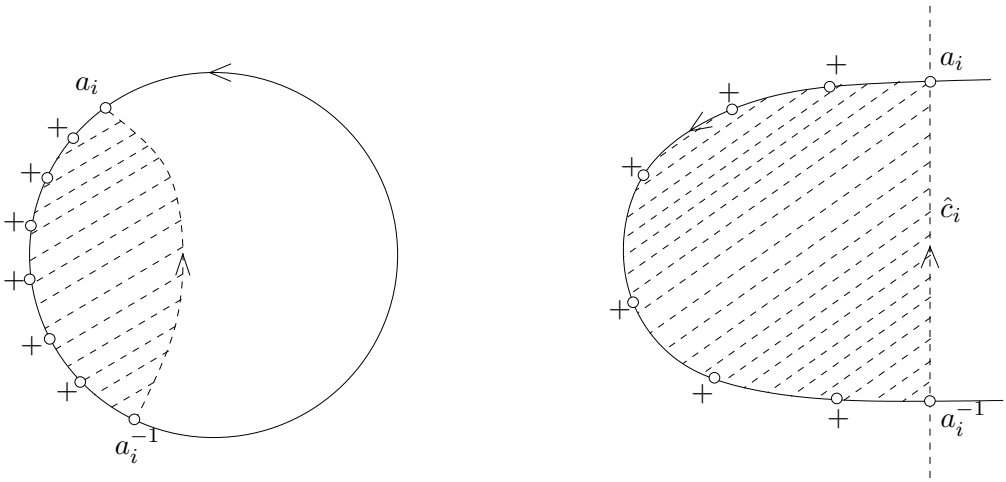}
\caption{Extending a pairing}
\label{fig:pairing_surround}
\end{figure}

\begin{lemma}
A pairing surrounds an embedded disc.
\label{lem:pairing}
\end{lemma}

\begin{figure}
\centering
\psfrag{cj}{{$\hat c_j$}}
\psfrag{ci}{{$\hat c_i$}}
\psfrag{pj}{{$p_j$}}
\psfrag{aj}{{$a_j$}}
\psfrag{aj-}{{$a_j^{-1}$}}
\psfrag{ai}{{$a_i$}}
\psfrag{ai-}{{$a_i^{-1}$}}
\includegraphics[width=.7\textwidth]{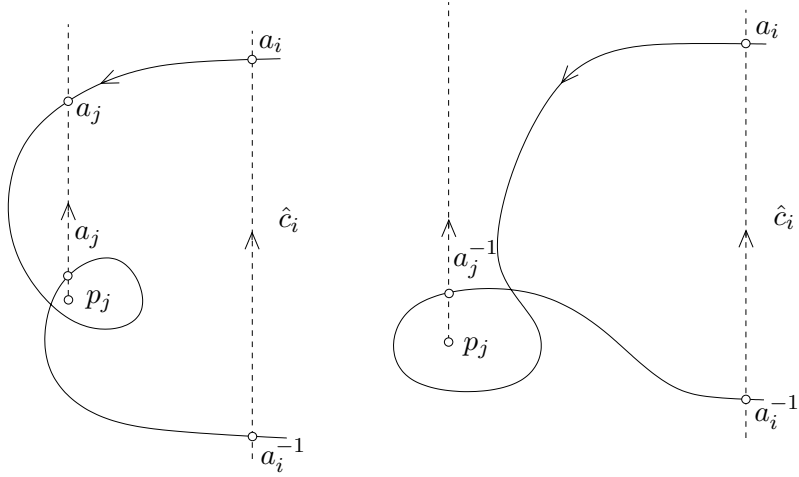}
\caption{The immersion inside of a pairing cannot have selfintersections.}
\label{fig:pairing}
\end{figure}

\begin{proof}
If the word $w(f)$ contains a pairing the immersion on $\S^2$ locally looks like the right side of Figure \ref{fig:pairing_surround} whereas the left side shows the situation in the preimage.\\
\\
Join the points $a_i$ and $a_i^{-1}$ in the preimage by an interval. The ray $\hat c_i$ in $\S^2$ induces an orientation on the interval such that it starts in $a_i^{-1}$ and ends in $a_i$. So we denote the interval by $[a_i^{-1},a_i]$. Since the rays $\hat c_j$ are chosen disjoint, the ray $\hat c_i$ has no further intersection points with $f(\S^1)$ or another ray between the points $a_i$ and $a_i^{-1}$. Therefore we can extend the immersion $f\colon \S^1 \to \S^2$ to an immersion $\widetilde f\colon \S^1 \cup [a_i^{-1},a_i] \to \S^2$.\\
Now take the restriction $\hat f\colon \S^1_{[a_i,a_i^{-1}]} \cup [a_i^{-1},a_i] \to \S^2$. We will show by contradiction that this restriction $\hat f$ is an embedding. Assume that it is not an embedding. Since the ray $\hat c_i$ is chosen with no selfintersections, $f(\S^1_{[a_i,a_i^{-1}]})$ need to have a selfintersection. Because $\S^1_{[a_i,a_i^{-1}]}$ is connected $f(\S^1_{[a_i,a_i^{-1}]})$ is connected as well. Hence there must be at least one loop in $f(\S^1_{[a_i,a_i^{-1}]})$. This loop encloses a component $X_j$, that is, a ray $\hat c_j$ starts at $p_j \in X_j$. This ray yields intersection points. 

We have to distinguish between two cases: The first case is that the loop  is on the left of $f(\S^1_{[a_i,a_i^{-1}]})$ and the second that it is on the right. In the first case the ray $\hat c_j$ intersects $f(\S^1_{[a_i,a_i^{-1}]})$ in two consecutive points $a_ja_j=a_j^2$ and in the second case it intersects $f(\S^1_{[a_i,a_i^{-1}]})$ in a negative point $a_j^{-1}$ (Figure \ref{fig:pairing}). Both times in contradiction to the assumption that there is a positive word inside of the pairing. Hence $\hat f$ is an embedding and extends to an embedded disc.
\end{proof}

\begin{lemma}
A negative group surrounds an embedded disc.
\label{lem:neg_group}
\end{lemma}

\begin{figure}
\centering
\psfrag{ci}{{$\hat c_i$}}
\psfrag{cj}{{$\hat c_j$}}
\psfrag{aj-}{{$a_j^{-1}$}}
\psfrag{ai-}{{$a_i^{-1}$}}
\psfrag{x0}{{$x_0$}}
\psfrag{pi}{{$p_i$}}
\psfrag{pj}{{$p_j$}}
\includegraphics[width=.9\textwidth]{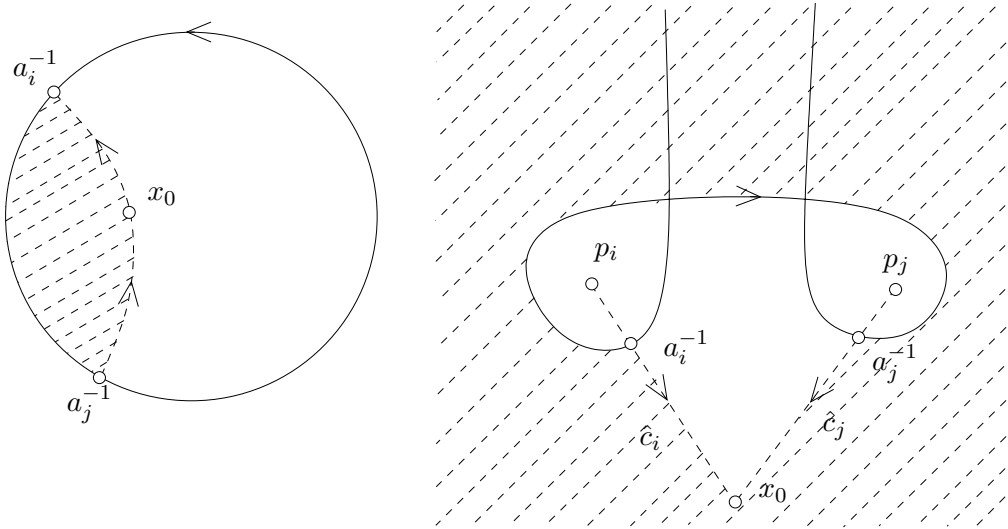}
\caption{Extending a negative group}
\label{fig:neg_group_surround}
\end{figure}

\begin{proof}
In the case of a negative group we locally get the situation depicted in Figure \ref{fig:neg_group_surround}.
Extend the immersion $f\colon \S^1 \to \S^2$ to the interval $[a_j^{-1},a_i^{-1}]$ and look at the restriction $\hat f\colon \S^1_{[a_i^{-1},a_j^{-1}]} \cup [a_j^{-1},a_i^{-1}] \to \S^2$. We will show that $\hat f$ is an embedding. As in the proof of Lemma \ref{lem:pairing} we assume that it is not, hence there has to be a double point in $f(\S^1_{[a_i^{-1},a_j^{-1}]} \cup [a_j^{-1},a_i^{-1}])$. Since $[a_j^{-1},a_i^{-1}]$ is mapped to $\hat c_j, \hat c_i$ and $x_0$ the double point has to be in $f(\S^1_{[a_i^{-1},a_j^{-1}]})$. Similar to the proof of Lemma \ref{lem:pairing} there has to be either consecutive intersection points $a_ka_k=a_k^2$ or a negative intersection point $a_k^{-1}$. %(Figure \ref{fig:neg_group}). 
But since no further intersection points occur in a negative group this is a contradiction and hence $\hat f$ is an embedding.
\end{proof}

Now we have seen that the special subwords, which were introduced in Definition \ref{def:subword}, surround embedded discs.

\begin{theorem}[Extension Theorem for Immersed Discs]
\index{Theorem!Extension Theorem!immersed discs}
Let $f\colon \S^1 \to \S^2$ be a normal immersion with word $w(f)$. If $\tau(f)=1-2\omega_1$ and $f$ has a grouping $\G$ then $f$ extends to an immersion $F:\overline D \to \S^2$.
\label{theo:extend}
\end{theorem}

\begin{figure}
\centering
\psfrag{a1}{{$a_1$}}
\psfrag{a2}{{$a_2$}}
\psfrag{a3}{{$a_3$}}
\psfrag{a1-}{{$a_1^{-1}$}}
\psfrag{a2-}{{$a_2^{-1}$}}
\psfrag{a3-}{{$a_3^{-1}$}}
\psfrag{Dj}{{\huge $D_j$}}
\psfrag{Dk}{{$D_k$}}

\includegraphics[width=.5\textwidth]{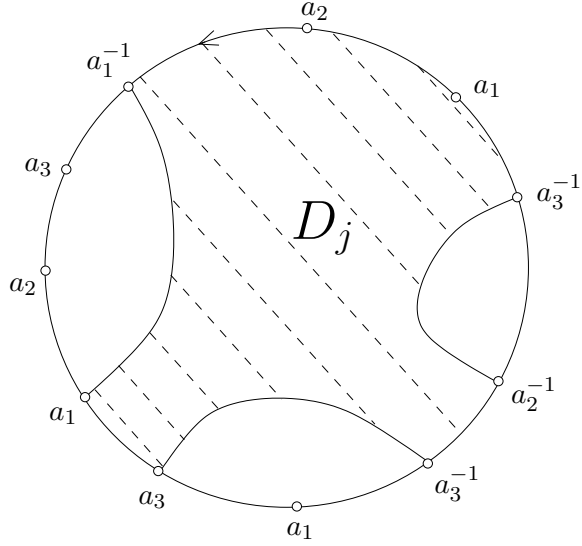}
\caption{Restriction of $f^*$ to a smaller disc $D_j$.}
\label{fig:restriction}
\end{figure}

\begin{proof}
Since $w(f)$ is groupable, a positive word remains after canceling $k_p$ pairings and $k_n$ negative groups, i.e., the grouping $\G$ has $k:=k_p+k_n$ edges (Lemma \ref{lem:tree}), each of them crossing an interval $I_1,\ldots,I_k$.
 
Denote the boundary points of the interval $I_j$ by $a_j$ and $b_j$, that is, $I_j=[a_j,b_j]$. These intervals decompose $\overline D$ in $k+1$ smaller discs $D_1,\ldots,D_{k+1}$. 

The boundary of the interval $I_j$ maps to the intersection points $a_j$ and $b_j$ of $f(\S^1)$ and a ray $\hat c_j$. W.l.o.g. assume that the ray $\hat c_j$ is oriented from $a_j$ to $b_j$. Denote that part of $\hat c_j$, which starts at $a_j$ and ends at $b_j$ with $\hat c_{[a_j,b_j]}=\hat c_{I_j}$. Then there are diffeomorphisms $\phi_j\colon I_j \to \hat c_{I_j}$, such that $f\colon \S^1 \to \S^2$ extends to $f^*\colon \S^1 \cup I_1 \cup \ldots \cup I_k \to \S^2$.\\
\\
Have a closer look at a small disc $D_j$. The boundary of $D_j$ contains $n$ Intervals $I_1,\ldots,I_n$. Since the intervals are disjoint and the boundary points of each interval are in $\S^1$, the boundary of $D_j$ contains $n$ connected subsets $S_1,\ldots,S_n \subset \S^1$ (Figure \ref{fig:restriction}). 

Each disc $D_j$ contains a word in its boundary. Since the intervals contain no letters this word is composed of the subwords of $w(f)$ which are contained in the subsets $S_k$. Because $w(f)$ is groupable the cancellation process leads to a positive word. Each negative letter of $w(f)$ belongs to a pairing or negative group and so is canceled. Thus the boundary of $D_j$ contains a positive word. 

Consider the restriction $f^{*}|_{\partial D_j} \colon \partial D_j \to \S^2$. Since $\partial D_j \simeq \S^1$ we can choose a system of rays for $f^{*}|_{\partial D_j}$. The rays $\hat c_j^*$ for $f^{*}|_{\partial D_j}$ are a selection of the rays $\hat c_j$ of $f$. Hence the resulting word $w(f^{*}|_{\partial D_j})$ is positive and
according to Lemma \ref{lem:positive_word} the restriction $f^{*}|_{\partial D_j}$ is an embedding which extends to an embedding with the help of the Riemann mapping Theorem (see \cite{Ahlfors}, p. 172ff). \\
\\
The Schwarz reflection principle assures that the restrictions $f^*$ can be glued together by an analytic transformation. That way the condition of the Sewing Theorem (see \cite{Courant}, Theorem 2.5) are fulfilled and hence the restrictions can be glued together to an immersion $F\colon\overline D \to \S^2$, which is the desired extension of $f\colon\S^1 \to \S^2$.
\end{proof}

\begin{prop}
Suppose $f\colon \S^1 \to \S^2$ is a normal immersion and $w(f)$ the corresponding reduced word. Denote the degree of the base point component by $\omega_1$. 
If $w(f)$ is groupable then each grouping
of $w(f)$ contains exactly $\omega_1$ negative groups.
\label{prop:preimagex0}
\end{prop}

\begin{proof}
Since $w(f)$ is groupable, the normal immersion $f\colon \S^1 \to \S^2$ extends to an immersion $F\colon \overline D \to \S^2$. We choose
the immersion which is constructed in Theorem \ref{theo:extend}.\\

Each negative group in the cancellation process of $w(f)$ induces an interval. Denote the intervals which are induced by negative groups by $I_1,\ldots,I_{n}$. 
According to the construction of $F$ for the preimage of $x_0$ we have
\[
F^{-1}(\{x_0\}) \subset \bigcup_{j=1}^n I_j.
\]
Since each interval contains exactly one preimage of $x_0$ (Lemma \ref{lem:neg_group}) and the degree of the base point component is $\omega_1$, it follows that $n=\omega_1$.
\end{proof}

Since $\omega_1$ only depends on $f\colon \S^1 \to \S^2$ (Proposition \ref{prop:preimage_const}) we get:

\begin{cor}
Every grouping of a reduced word $w(f)$ has the same number of negative groups.
\label{cor:neggroup_const}
\end{cor}

\begin{bem}
\emph{Samuel J. Blank \cite{Blank} analyzed normal immersions $f\colon\S^1 \to \R^2$. In that case there is an outer component with degree $0$ for every normal immersion, since an extension 
has to be compact and $\R^2$ is not compact. According to Proposition \ref{prop:preimagex0} no negative groups occur in the grouping. Hence only pairings are canceled out, which is
in fact the same as what Samuel J. Blank showed in his proof of his theorem in the case of normal immersions $f\colon\S^1 \to \R^2$.
So the proof of Theorem \ref{theo:extend} includes as a special case that each normal immersion $f\colon \S^1 \to \R^2$ can be written as a normal immersion $\widetilde f\colon\S^1 \to \S^2$ via 
stereographic projection.}
\end{bem}

\subsection{Ungroupable Words}

We have shown in the previous subsection that a normal immersion $f\colon\S^1 \to \S^2$ extends to an immersion $F\colon\overline D \to \S^2$ when the corresponding reduced word $w(f)$ is groupable. Now we will have a look at which properties of the reduced word $w(f)$ ban such an extension.

\begin{lemma}
Let $f\colon\S^1 \to \S^2$ be a normal immersion and $w(f)$ the corresponding reduced word. If $\widetilde{w}(f)$ is the word which
remains after canceling all pairings and negative groups then $\widetilde{w}(f)$ does not contain a subword $v=a_j^n$ with $n \geq 3$.
\label{lem:aj3}
\end{lemma}

\begin{figure}
\centering
\psfrag{pj}{{$p_j$}}
\psfrag{pk}{{$p_k$}}
\psfrag{aj1}{{$a_j^{(i)}$}}
\psfrag{aj2}{{$a_j^{(i+1)}$}}
\psfrag{aj3}{{$a_j^{(i+2)}$}}
\psfrag{ak}{{$a_k$}}
\psfrag{d1}{{$d_1$}}
\psfrag{d2}{{$d_2$}}
\psfrag{l}{{$\lambda$}}

\subfigure[Situation of consecutive intersection points]{\includegraphics[width=.45\textwidth]{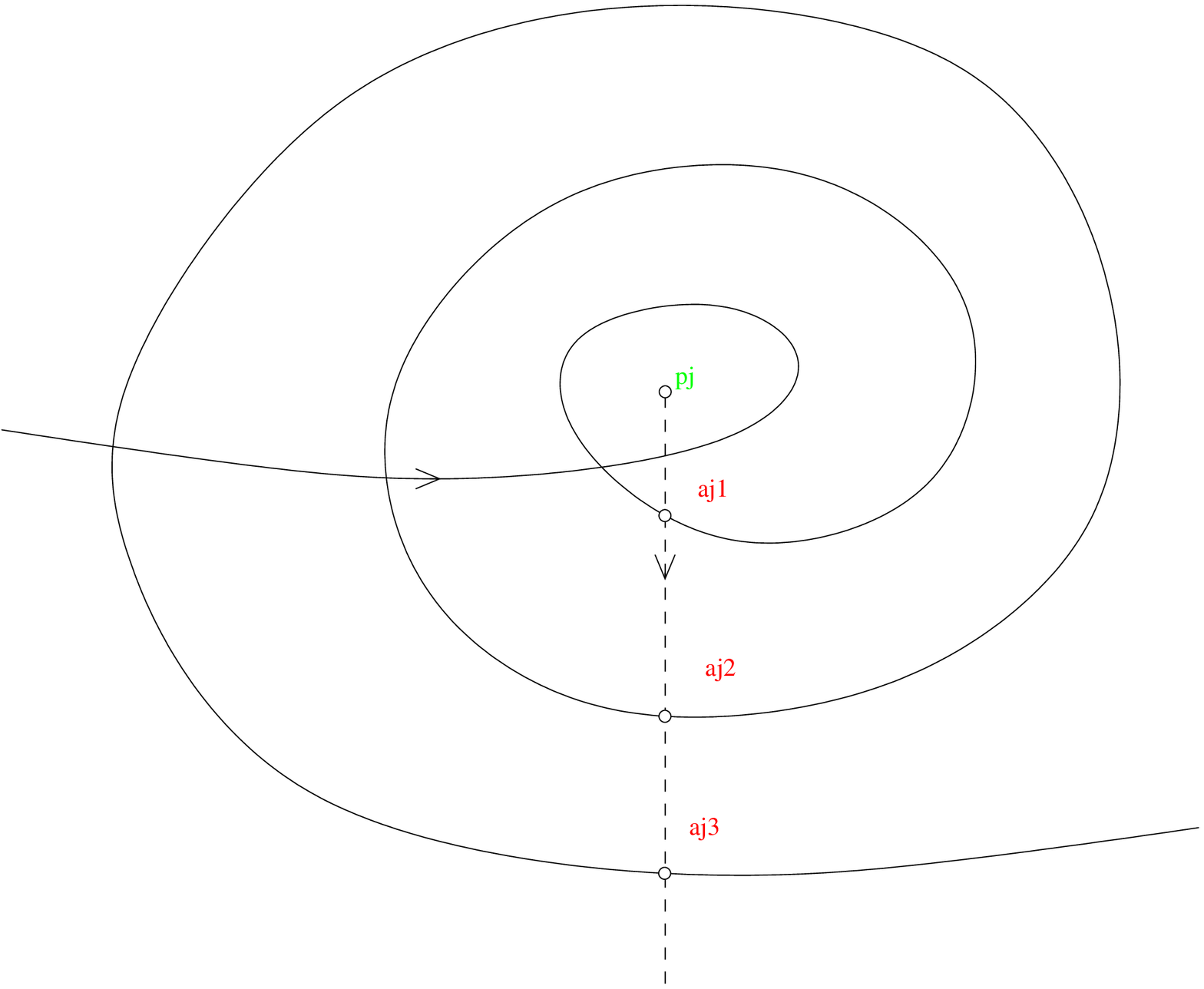}}\quad
\subfigure[The closed path $\lambda$ leads to another component]{\includegraphics[width=.45\textwidth]{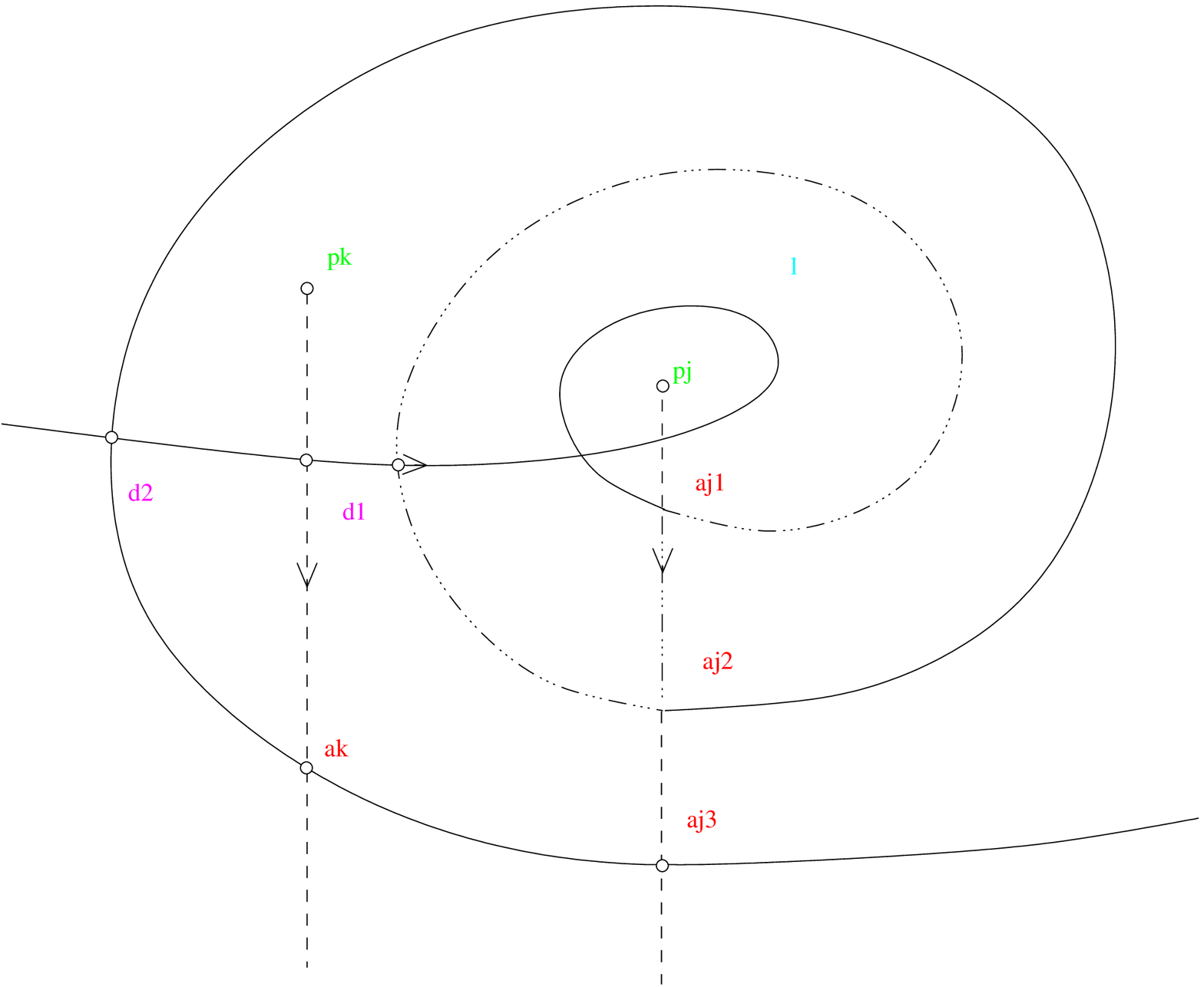}}
\caption{Consecutive Intersection points}
\label{fig:lambda}
\end{figure}

\begin{proof}
To get consecutive intersection points $a_j^n$ the immersion $f(\S^1)$ has to orbit the point $p_j$, because after the first intersection point $a_j$ the image $f(\S^1)$ has to intersect the ray $\hat c_j$ again in the same direction without intersecting another ray $\hat c_i$ (Figure \ref{fig:lambda} (a)). 

We denote the intersection points by $a_j^{(1)},\ldots,a_j^{(n)}$, according to their appearance in $\omega$. Two consecutive intersection points $a_{j}^{(i)}$ and $a_{j}^{(i+1)}$ define a closed curve $\lambda$, which starts at $a_{j}^{(i)}$, travels along $f(\S^1)$ to $a_{j}^{(i+1)}$ and along $\hat c_j$ back to $a_{j}^{(i)}$ (Figure \ref{fig:lambda} (b)). 
The first intersection point $a_j^{(1)}$ arises when the immersion $f(\S^1)$ intersects the ray $\hat c_j$ for the first time. Therefore this first intersection point lies on the same side of $\lambda$ as $p_j$. Since $f(\S^1)$ is connected it has to intersect $\lambda$ at another point. This could be on the part of $\lambda$ which belongs to $f(\S^1)$ or on the part which belongs to $\hat c_j$.

If this intersection point arises on the part which belongs to $\hat c_j$ then it produces an intersection point $a_j^{\pm 1}$. If it has the same sign as the one in $\omega$, then we have found another intersection point before $a_j^{(1)}$, which contradicts the choice of $a_j^{(1)}$ as the first point. Hence it has to be the other sign. But in that case we have a sequence $a_ja_j^{-1}$, which is a contradiction to the fact that the word is reduced. Therefore the intersection point has to arise on the part of $\lambda$ which belongs to $f(\S^1)$ and so it is a double point of $f(\S^1)$.

Since we can choose consecutive intersection points in $n-1$ ways, we get $n-1$ closed curves $\lambda_1,\ldots,\lambda_{n-1}$, each producing a double point of $f(\S^1)$. All double points have to lie on one side of the ray $\hat c_j$, otherwise there would be a new intersection point $a_j^{\pm 1}$. That means between every two consecutive intersection points there is a double point. Start at an intersection point $a_j^{(i)}$ and travel along $f(\S^1)$ until a double point $d_1$ is reached. Follow $f(\S^1)$ along the branch that passes through the next intersection point $a_j^{(i+1)}$ until we reach the next double point $d_2$. Since there is another branch of $f(\S^1)$ which joins the double points without passing $a_j^{(i+1)}$, this curve surrounds a component $X_k$ of $\S^2 \backslash f(\S^1)$ (Figure \ref{fig:lambda} (b)). Hence there must be a ray $\hat c_k$ which starts in that component. Since the ray $\hat c_k$ cannot intersect the ray $\hat c_j$ it has to intersect $f(\S^1)$ somewhere between the double point $d_2$ and the next intersection point $a_j^{(i+2)}$. That is, a new intersection point $a_k$ arises before $a_j^{(i+2)}$.

Hence there cannot be more then two equal consecutive intersection points, i.e., $v=a_j^n$ can only occur in $w(f)$ with $n\leq 2$. 
\end{proof}

\begin{figure}
\centering
\psfrag{x0}{{$x_0$}}
\psfrag{p1}{{$p_1$}}
\psfrag{p2}{{$p_2$}}
\psfrag{p3}{{$p_3$}}
\psfrag{p4}{{$p_4$}}
\psfrag{p5}{{$p_5$}}
\psfrag{p6}{{$p_6$}}
\psfrag{p7}{{$p_7$}}
\psfrag{a1}{{$a_1$}}
\psfrag{a2}{{$a_2$}}
\psfrag{a3}{{$a_3$}}
\psfrag{a4}{{$a_4$}}
\psfrag{a5}{{$a_5$}}
\psfrag{a6}{{$a_6$}}
\psfrag{a7}{{$a_7$}}
\psfrag{a2-}{{$a_2^{-1}$}}
\psfrag{a6-}{{$a_6^{-1}$}}
\psfrag{a7-}{{$a_7^{-1}$}}
\psfrag{D}{{\huge $\widetilde D$}}
\subfigure[Decomposing $\overline D$]{\includegraphics[width=.67\textwidth]{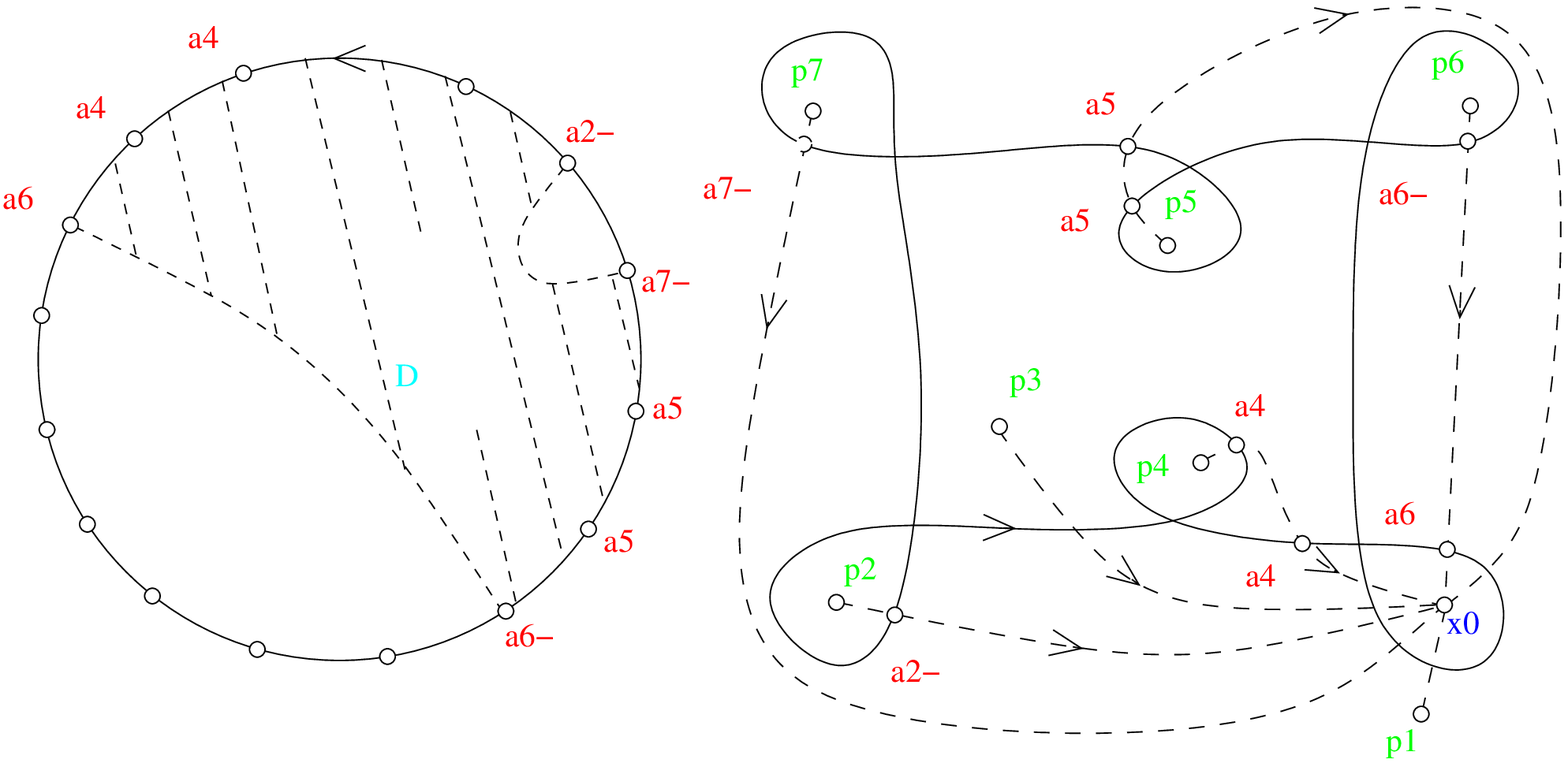}}\quad
\subfigure[Situation for $\widetilde D$]{\includegraphics[width=.3\textwidth]{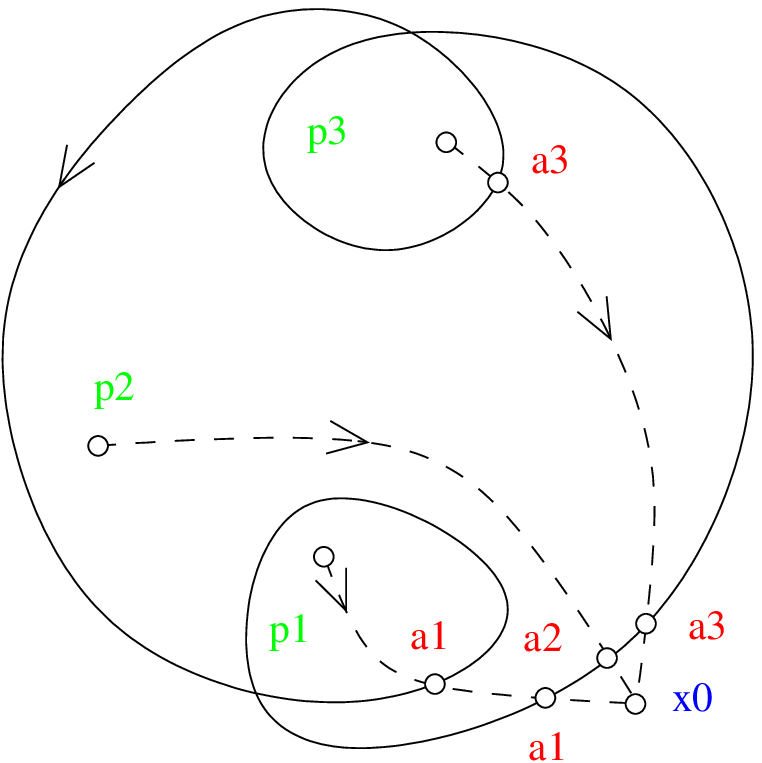}}
\caption{The remaining disc $\widetilde D$, containing subwords $a_j^2$.}
\label{fig:combination}
\end{figure}

\begin{lemma}
Let $f\colon\S^1 \to \S^2$ be a normal immersion with word $w(f)$. Denote by $\widetilde{w}(f)$ the word which remains after canceling all pairings and negative groups.
If $\widetilde{w}(f)$ contains a subword $v=a_j^2$ then $f$ cannot be extended.
\label{lem:aj2}
\end{lemma}

\begin{proof}
Take the reduced word $w(f)$ and cancel all pairings and negative groups. Each time a pairing or a negative group is canceled an interval is inserted which bounds a small disc together with a part of $\S^1$. According to Lemmas \ref{lem:pairing} and \ref{lem:neg_group} this surrounds an immersed disc. After canceling all pairings and negative groups a disc $\widetilde D$ remains, which has $\widetilde w(f)$ as a word in the boundary (Figure \ref{fig:combination}).

Assume that there are $n$ subwords $a_1^2,\ldots,a_n^2$ in $\widetilde{w}(f)$. Then the restriction $\widetilde f$ of $f$ to $\widetilde D$ is a circle with $n$ loops added (Figure \ref{fig:combination} (b)).
Decomposition of the immersion leads to $n+1$ embeddings. W.l.o.g. embedding $\widetilde f_0$ is the circle and embeddings $\widetilde f_1,\ldots,\widetilde f_n$ are the added loops.
Then the following is true for the $\widetilde f$
\[
\tau(\widetilde f)=\sum_{j=0}^n \tau(\widetilde f_j)=\tau(f_0)+\sum_{j=1}^n \tau(\widetilde f_j)=1+\sum_{j=1}^n \tau(\widetilde f_j).
\]
Since each loop produces two consecutive positive letters, the tangent winding number of each embedding $\widetilde f_1,\ldots,\widetilde f_n$ is $\tau(f_j)=1$. 
Hence 
\[
\tau(\widetilde f)=1+\sum_{j=1}^n \tau(\widetilde f_j)=1+n,
\]
and according to Proposition \ref{prop:notwendig} this immersion cannot be extended for $n \neq 0$.
\end{proof}

\begin{lemma}
Let $f\colon\S^1 \to \S^2$ be a normal immersion with word $w(f)$. Denote by $\widetilde{w}(f)$ the word which
remains after canceling all pairings and negative groups. If $\widetilde{w}(f)$ contains negative letters but no subword $v=a_j^2$ then $f$ cannot be extended.
\label{lem:aj-}
\end{lemma}

\begin{figure}
\centering
\psfrag{aj-}{{$a_j^{-1}$}}
\psfrag{al-}{{$a_l^{-1}$}}
\psfrag{ak-}{{$a_k^{-1}$}}
\psfrag{am}{{$a_m$}}
\psfrag{an}{{$a_n$}}
\psfrag{ai}{{$a_i$}}
\psfrag{pj}{{$p_j$}}
\psfrag{pl}{{$p_l$}}
\psfrag{pk}{{$p_k$}}
\psfrag{D}{{\huge $\widetilde D$}}
\psfrag{f}{$f|_{\partial \widetilde D}$}
\includegraphics[width=\textwidth]{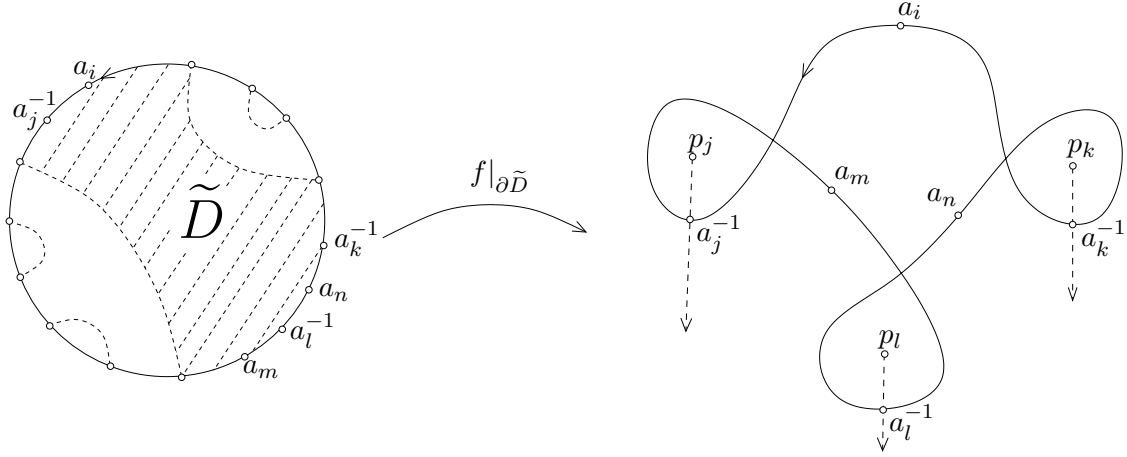}
\caption{The remaining disc $\widetilde D$ with single loops added}
\label{fig:aj-}
\end{figure}

\begin{proof}
Assume $\widetilde{w}(f)$ contains $n$ negative letters. Since $\widetilde{w}(f)$ contains no negative groups, each negative letter is separated by at least one positive letter. 

Canceling all pairings and negative groups leads to a decomposition of $\overline D$ and a remaining disc $\widetilde D$ with $\widetilde{w}(f)$ in the boundary. Since all negative letters are separated and no subword $a_j^2$ occurs the restriction $f|_{\partial \widetilde D}$ maps to a circle with $n$ single loops added (Figure \ref{fig:aj-}). Decomposing $f|_{\partial \widetilde D}$ into embeddings  leads to $n+1$ embeddings $f_0,\ldots,f_n$. W.l.o.g. $f_0$ is mapped to the circle and $f_1,\ldots,f_n$ to the added loops. Since $f_1,\ldots,f_n$ produce negative letters the tangent winding number is $\tau(f_j)=-1$, for $j=1,\ldots,n$. Thus for the tangent winding number of $f|_{\partial \widetilde D}$ we have
\[
\tau(f|_{\partial \widetilde D})=\sum_{j=0}^n \tau(f_j)=\tau(f_0) + \sum_{j=1}^n \tau(f_j) \overset{\tau(f_j)=-1}=\tau(f_0)-n.
\]
Since $f_0$ maps to the circle $\tau(f_0)=+1$ and therefore $\tau(f|_{\partial \widetilde D})=1-n$. According to Proposition \ref{prop:notwendig} $f$ cannot be extended if $n$ is odd.\\
\\
If $n$ is even then Proposition \ref{prop:notwendig} shows that the degree $\omega_1$ of the base point component is
\[
\omega_1={n \over 2}.
\]
If $f|_{\partial \widetilde D}$ can be extended then $\widetilde{w}(f)$ is groupable and hence contains $\omega_1$ negative groups \linebreak (Proposition \ref{prop:preimagex0}). Since $\widetilde{w}(f)$ contains no negative groups by assumption, $f|_{\partial \widetilde D}$ fails to extend and hence $f$ cannot be extended.
\end{proof}

Now we can show that a normal immersion $f\colon \S^1 \to \S^2$ with an ungroupable reduced word $w(f)$ cannot be extended to an immersion $F\colon \S^1 \to \S^2$.

\begin{theorem}
Let $f\colon \S^1 \to \S^2$ be a normal immersion and $w(f)$ the corresponding reduced word. 
If $\tau(f) \neq 1-2\omega_1$ or $w(f)$ is not groupable then $f$ cannot be extended to an immersion $F\colon\overline D \to \S^2$.
\label{theo:ungroupable}
\end{theorem}

\begin{proof}
If the tangent winding number $\tau(f)$ of $f$ is not equal to $1-2\omega_1$ then Proposition \ref{prop:notwendig} shows that $f$ cannot be extended.

So assume $f$ has the tangent winding number $\tau(f)=1-2\omega_1$ but the word $w(f)$ is not groupable. Cancel all pairings and negative groups of $w(f)$ and denote the remaining word by $\widetilde{w}(f)$. Since $w(f)$ is not groupable
$\widetilde{w}(f)$ is not a positive word. 

If $\widetilde{w}(f)$ contains a subword $a_j^2$ then $f$ cannot be extended according to Lemma \ref{lem:aj2}. Hence we can assume that $\widetilde{w}(f)$ contains no subword $a_j^2$. But since $\widetilde{w}(f)$ is not a positive word it has to contain negative letters. In this case $f$ cannot be extended according to Lemma \ref{lem:aj-}.
\end{proof}

\subsection{Uniqueness of Extensions}
\label{subsec:different_ext}

Recently we have seen that a normal immersion $f\colon \S^1 \to \S^2$ can be extended to an immersion $F\colon \overline D \to \S^2$ if and only if the corresponding reduced word $w(f)$ is groupable. In this section we will show that the number of different extensions is equivalent to the number of different groupings.

\begin{defs}
\emph{Suppose $f\colon \S^1 \to \S^2$ is a normal immersion and  $F_1, F_2\colon \overline D \to \S^2$ are extensions of $f$. Two extensions are \emph{equivalent}\index{extension!equivalent} if there exists an orientation preserving diffeomorphism $\phi\colon \overline D \to \overline D$ such that $F_1=F_2 \circ \phi$ and $\phi|_{\S^1}=\id_{\S^1}$.}
\label{def:ext_equal}
\end{defs}

\begin{bem}
\emph{In the proof of Theorem \ref{theo:extend} an extension is constructed to a given normal immersion with groupable reduced word $w(f)$. In fact, this extension defines an equivalence class of extensions. If we talk about extensions then we talk about the equivalence classes. As representative we pick an extension which is constructed in the proof of Theorem \ref{theo:extend}.}
\end{bem}

\begin{theorem}[Uniqueness Theorem for Immersed Discs]
\index{Theorem!Uniqueness Theorem!immersed discs}
Suppose $f\colon \S^1 \to \S^2$ is a normal immersion and $w(f)$ the corresponding reduced word. Then two different groupings define two different equivalence classes of extensions of $f\colon\S^1 \to \S^2$. 
\label{theo:different}
\end{theorem}

\begin{proof}
Theorem \ref{theo:extend} shows that a grouping $\G$ defines an equivalence class of extensions. Assume $w(f)$ has two different groupings $\G_1, \G_2$, i.e., the weighted trees are not isomorphic. \\
\\
We will show the claim by contradiction. Thus assume that $F_1,F_2\colon \S^1 \to \S^2$ are two equivalent extensions of $f\colon\S^1 \to \S^2$. Then an orientation preserving diffeomorphism $\phi\colon\overline D \to \overline D$ exists, such that $F_1=F_2 \circ \phi$ and $\phi|_{\S^1}=\id_{\S^1}$.

Since $f$ extends the reduced word $w(f)$ has a grouping $\G_1$ which induces a decomposition of $\overline D$ by intervals $\mathcal{I}_1$. Assume that this grouping induces the extension $F_1$. Since $\phi|_{\S^1}=\id_{\S^1}$ the intervals $\mathcal{I}_1$ are mapped to homotopic intervals $\mathcal{I}_2$ by $\phi$. Thus the decomposition of $\overline D$ by intervals $\mathcal{I}_2$ induce the extension $F_2$. Since the intervals $\mathcal{I}_2$ are still disjoint, these intervals $\mathcal{I}_2$ induce a tree $\G_2$ which is isomorphic to the tree of $\G_1$. Since $\phi$ does not change the boundary of the intervals the weighted trees are isomorphic as well.\\ 
\\
Hence different weighted trees, and thus different groupings, define nonequivalent extensions.
\end{proof}

\bibliography{biblio}{}

\begin{thebibliography}{GKS07}

\bibitem[Ahl53]{Ahlfors}
Lars~V. Ahlfors.
\newblock {\em {Complex analysis, an introduction to the theory of analytic
  functions of one complex variable.}}
\newblock {International Series in Pure and Applied Mathematics. London:
  McGraw-Hill Publishing Co., Ltd. XI, 247 p. }, 1953.

\bibitem[Bla67]{Blank}
Samuel~Joel Blank.
\newblock {\em Extending Immersions and regular Homotopies in Codimension $1$}.
\newblock PhD thesis, Brandeis University, 1967.

\bibitem[Cou50]{Courant}
Richard Courant.
\newblock {\em {Dirichlet's principle, conformal mapping and minimal
  surfaces.}}
\newblock {(Pure and Applied Mathematics. A Series of Texts and Monographs.
  III) New York: Interscience Publishers, Inc. XIII, 330 p. }, 1950.

\bibitem[Fra73]{Francis2}
George~K. Francis.
\newblock {Spherical curves that bound immersed discs.}
\newblock {\em Proc. Am. Math. Soc.}, 41:87--93, 1973.

\bibitem[Fri10]{Frisch}
Dennis Frisch.
\newblock {\em {Classification of Immersions which are Bounded by Curves in
  Surfaces}}.
\newblock PhD thesis, TU Darmstadt, 2010.
\newblock {http://tuprints.ulb.tu-darmstadt.de/2202/}.

\bibitem[GKS07]{GKS}
Karsten {Grosse-Brauckmann}, Robert~B. Kusner, and John~M. Sullivan.
\newblock {Coplanar constant mean curvature surfaces.}
\newblock {\em Commun. Anal. Geom.}, 15(5):985--1023, 2007.

\bibitem[MC93]{MC}
Margaret McIntyre and Grant Cairns.
\newblock {A new formula for winding number.}
\newblock {\em Geom. Dedicata}, 46(2):149--159, 1993.

\bibitem[Po{\'e}69]{Poenaru}
Valentin Po{\'e}naru.
\newblock {Extension des immersions en codimension 1 (d'apr\`es Samuel Blank).
  (Extension of immersions in codimension 1 (following Samuel Blank)).}
\newblock {Sem. Bourbaki 1967/1968, Expose No.342, 33 p. (1969).}, 1969.

\end{thebibliography}
\bibliographystyle{alpha}

\end{document}